\documentclass[11pt,reqno]{amsart}

\usepackage{amsmath,amssymb,amsthm}
\usepackage{bm}
\usepackage{mathrsfs}
\usepackage[all]{xy}
\usepackage{booktabs}
\usepackage{fullpage}
\usepackage[pagebackref=true]{hyperref}
\usepackage{mathtools}
\usepackage[abbrev,msc-links]{amsrefs}
\usepackage{caption,subcaption}
\usepackage[]{todonotes}
\hypersetup{
 hidelinks,
 bookmarksopen=true,
}

\newcommand{\itg}{\mathbb{Z}}

\newcommand{\rl}{\mathbb{R}}

\newcommand{\CB}{\mathcal{B}}

\newcommand{\CF}{\mathcal{F}}

\newcommand{\CS}{\mathcal{S}}
\newcommand{\CT}{\mathcal{T}}

\newcommand{\Fg}{\mathfrak{g}}
\newcommand{\Fh}{\mathfrak{h}}

\newcommand{\Fk}{\mathfrak{k}}
\newcommand{\Fl}{\mathfrak{l}}

\newcommand{\Fn}{\mathfrak{n}}
\newcommand{\Fo}{\mathfrak{o}}
\newcommand{\Fp}{\mathfrak{p}}

\newcommand{\Fv}{\mathfrak{v}}

\newcommand{\Fz}{\mathfrak{z}}

\newcommand{\isom}{\xrightarrow{\sim}}

\theoremstyle{plain}
\newtheorem{theorem}{Theorem}[section]

\newtheorem{lemma}[theorem]{Lemma}
\newtheorem{proposition}[theorem]{Proposition}
\newtheorem{corollary}[theorem]{Corollary}
\theoremstyle{definition}
\newtheorem{definition}[theorem]{Definition}
\newtheorem{conjecture}[theorem]{Conjecture}

\theoremstyle{definition}
\newtheorem{remark}[theorem]{Remark}

\begin{document}

\title{A Hilbert--Mumford criterion for nilsolitons}
\author[Y.~Hashimoto]{Yoshinori Hashimoto}
\date{\today}
\address{Department of Mathematics, Osaka Metropolitan University, 3-3-138, Sugimoto, Sumiyoshi-ku, Osaka, 558-8585, Japan.}
\email{yhashimoto@omu.ac.jp}

\begin{abstract}
We give an algebraic criterion for a nilpotent real Lie algebra and prove that it provides a necessary and sufficient condition for the associated nilpotent Lie group to admit left-invariant Ricci solitons, called nilsolitons. As an application of this result, we generalise Nikolayevsky's criterion for the existence of nilsolitons to nilpotent Lie algebras without nice bases. We further prove a modified version of the Taketomi--Tamaru conjecture for nilpotent Lie groups which gives an obstruction to the existence of nilsolitons.
\end{abstract}

\maketitle

\tableofcontents

\section{Introduction}

\subsection{Statement of the main results}

Finding a canonical Riemannian metric on a given manifold is a very important topic in Riemannian geometry. An intensively studied problem is whether there exists an Einstein metric, or more generally a Ricci soliton. Such a metric may or may not exist on a given Riemannian manifold, making it a very interesting and difficult problem.

In this paper we focus on the case when the manifold is a homogeneous nilmanifold, which is a nilpotent Lie group $N$ equipped with a left-invariant metric; recall that a connected and simply connected Lie group $N$ is said to be nilpotent if its Lie algebra $\Fn$ is a nilpotent Lie algebra over $\rl$. We assume throughout that $\Fn$ is not commutative. A classical result (see \cite[Theorem 4]{Jen69} and \cite[Theorem 2.4]{Mil76}) shows that $N$ cannot admit a left-invariant Einstein metric. We thus consider left-invariant Ricci solitons, also called \textit{nilsolitons}, which may or may not exist on $N$. The existence and non-existence of nilsolitons, and their classification, have been studied intensively for many years. Furthermore, the resolution of the Alekseevskii conjecture \cite{BL23} means that the classification of non-compact homogeneous Einstein spaces reduces to that of nilsolitons (see also \cite{BL22,Lau09,Heb98}).

The main result of this paper is a numerical algebraic criterion which gives a necessary and sufficient condition for the existence of nilsolitons. The criterion is given by the following invariant which is defined for an orbit of a one-parameter subgroup emanating from the Lie bracket $\mu \in \wedge^2 \Fn^* \otimes \Fn$ on $\Fn$, analogously to the Hilbert--Mumford weight in Geometric Invariant Theory (GIT) \cite{MFK94,Kirw}. We first fix an inner product $\langle , \rangle$ on $\Fn$ so that a distinguished derivation $\phi_{\mu}$ of $(\Fn , \mu)$, called the \textit{pre-Einstein derivation} \cite{Nik11}, is self-adjoint with respect to it; see \S \ref{scsln} for more details. For an endomorphism $\lambda \in \mathrm{End} (\Fn)$ which is self-adjoint with respect to $\langle , \rangle$, we define the Hilbert--Mumford weight $\nu (\lambda ; \mu )$ as
	\begin{equation*}
		\nu(\lambda  ; \mu ) := - \min \{ \lambda_i+\lambda_j - \lambda_k  \mid 1 \le i < j \le n, 1 \le k \le n, \text{ s.t.~} \mu^k_{ij} \neq 0 \},
	\end{equation*}
	where we presented $\lambda$ as a diagonal matrix $\mathrm{diag}(\lambda_1 , \dots , \lambda_n)$ by choosing a $\langle, \rangle$-orthonormal basis $\{ e_i \}_{i=1}^n$ for $\Fn$ with respect to which we write $\mu (e_i, e_j) = \sum_{k=1}^n \mu^k_{ij} e_k$. We compute this invariant relative to the subspace defined by
	\begin{equation*}
		\CS := \rl (I_{\Fn} - \phi_{\mu}) \subset \mathrm{End} (\Fn ), 
	\end{equation*}
	in the sense that we only consider $\lambda \in \mathrm{End} (\Fn)$ that commutes with $\CS$ and is orthogonal to $\CS$, where $I_{\Fn}$ stands for the identity endomorphism on $\Fn$ and we use the inner product $\langle , \rangle_{\mathrm{End} (\Fn)}$ on $\mathrm{End} (\Fn)$ naturally induced from $\langle , \rangle$ to define the orthogonality. The space $\CS$ is proportional to the Ricci endomorphism of the nilsoliton if there exists one (Theorem \ref{thnikrecpe}).
	
	The precise statement of our main result is as follows.

\begin{theorem} \label{thmain}
	Let $(N, \textsl{g}_N)$ be a homogeneous nilmanifold and $(\Fn , \mu )$ be the corresponding nilpotent Lie algebra with the inner product $\langle , \rangle$, with respect to which the pre-Einstein derivation is self-adjoint. Then, $N$ admits a left-invariant Ricci soliton if and only if $\nu(\lambda  ; \mu )  \ge 0$ for any self-adjoint endomorphism $\lambda \in \mathrm{End} (\Fn)$ which commutes with $\CS$ and is orthogonal to $\CS$ with respect to $\langle , \rangle_{\mathrm{End} (\Fn)}$, with equality if and only if $\lambda$ is a derivation of $\mu$.
\end{theorem}

The result above can be regarded as a nilsoliton version of the Kempf--Ness theorem \cite{KemNes78}, which establishes the link between zero of moment maps and GIT stability for complex smooth projective varieties. The proof that we give in this paper is inspired by the original complex setting, which uses the variational principle involving a geodesically convex energy functional, but the situation is closer to the setting of the relative stability in GIT that we briefly review in \S \ref{scrtcs}; see also Remark \ref{rmrsfrs} for more details. For all these arguments to work, it is crucial that the Ricci curvature of a homogeneous nilmanifold can be regarded as a moment map, as pointed out in the foundational works of Lauret \cite{Lau01,Lau01b,Lau03,Lau06,Lau09}.

A result similar but inequivalent to Theorem \ref{thmain} was already proved by Jablonski \cite[Section 6 and Theorem 6.1]{Jab11}. While the differences are explained in detail in \S \ref{scrtcs}, we present how Theorem \ref{thmain} can be applied to prove new results, with the hope of clarifying its novelty, rather than making detailed comparisons here. We first prove the following generalisation of theorems by Nikolayevsky (\cite[Theorem 3]{Nik11} and \cite[Theorem 1 and Corollary 1]{Nik08ded}, see also Remark \ref{rmnkscal}) that are known to be helpful in dealing with concrete examples.

\begin{theorem} \label{thnik}
	Let $(\Fn, \mu , \langle , \rangle)$ be a metric nilpotent Lie algebra such that the pre-Einstein derivation $\phi_{\mu}$ is self-adjoint with respect to $\langle , \rangle$, and let $\CB := \{ e_i \}_{i=1}^n$ be a $\langle ,\rangle$-orthonormal basis for $\Fn$ which diagonalises $\phi_{\mu}$. We define a finite subset $\CF_{\CB}$ of $\Fn$ by
	\begin{equation*}
		\CF_{\CB} := \{ e_i + e_j - e_k \mid 1 \le i < j \le n, 1 \le k \le n, \textup{ s.t.~} \mu^k_{ij} \neq 0 \} ,
	\end{equation*}
	where we wrote $\mu(e_i, e_j) = \sum_{k=1}^n \mu_{ij}^k e_k$ as before, and set $m := |\CF_{\CB} |$ to label elements of $\CF_{\CB}$ as $f_1 , \dots , f_m$. Let $\mathsf{L}_{\CB}$ be the affine span of $\CF_{\CB}$ (i.e.~the smallest affine subspace of $\Fn$ containing $\CF_{\CB}$), $\mathrm{Conv} ( \CF_{\CB})$ be the convex hull of $\CF_{\CB}$, and $P_0$ be the projection of the origin of $\Fn$ to $\mathsf{L}_{\CB}$, whose length with respect to $\langle , \rangle$ turns out to be non-zero and is denoted by $\Vert P_0 \Vert$. Then the following hold.
	\begin{enumerate}
		\item Let $\Fp_{\CB}$ be the set consisting of endomorphisms on $\Fn$ that are diagonal with respect to $\CB$ and orthogonal to $\CS$ with respect to $\langle , \rangle_{\mathrm{End} (\Fn)}$. Then, $\nu (\lambda ; \mu ) > 0$ for any $\lambda \in \Fp_{\CB} \setminus \mathrm{Der} (\mu )$ if and only if $P_0$ lies in the interior of $\mathrm{Conv} (\CF_{\CB})$.
		\item Define an $m \times n$ matrix $\mathbf{Y}_{\CB}$ whose $(p,q)$-th entry is given by the following: when $f_p \in \CF_{\CB}$ is written as $f_p= e_i + e_j - e_k$, we set
		\begin{equation*}
			(\mathbf{Y}_{\CB})_{pq} := \begin{cases}
				1 &\quad (q=i, j) \\
				-1 &\quad (q=k) \\
				0 &\quad (\textup{otherwise})
			\end{cases}.
		\end{equation*}
		Then, $P_0$ lies in the interior of $\mathrm{Conv} (\CF_{\CB})$ if and only if there exist $\alpha_1 , \dots , \alpha_m >0$ such that $\sum_{i=1}^m \alpha_i = \Vert P_0 \Vert^{-2}$ and
		\begin{equation} \label{eqnikyyt}
			\mathbf{Y}_{\CB} \mathbf{Y}_{\CB}^t \begin{pmatrix}
				\alpha_1 \\
				\vdots \\
				\alpha_m
			\end{pmatrix}
			=\begin{pmatrix}
				1 \\
				\vdots \\
				1
			\end{pmatrix}.
		\end{equation}
	\end{enumerate}
	In particular, given Theorem \ref{thmain}, a homogeneous nilmanifold $(N , \textsl{g}_{N})$ associated to $(\Fn, \mu , \langle , \rangle)$ admits a left-invariant Ricci soliton if and only if $P_0$ lies in the interior of $\mathrm{Conv} (\CF_{\CB})$, or equivalently (\ref{eqnikyyt}) admits a solution $\alpha_1 , \dots , \alpha_m >0$ with $\sum_{i=1}^m \alpha_i = \Vert P_0 \Vert^{-2}$, for any $\langle , \rangle$-orthonormal basis $\CB$ which diagonalises $\phi_{\mu}$.
\end{theorem}

Nikolayevsky proved the result above when $(\Fn , \mu )$ has a \textit{nice} basis \cite[Definition 3]{Nik11}, for which the statement is simpler as it suffices to check the conditions for just one nice basis (see e.g.~\cite[Proof of Lemma 2]{Nik11}). We need to check the conditions in (i) for all $\langle , \rangle$-orthonormal bases for it to be generalised to nilpotent Lie algebras without nice bases, such as \cite[Examples 3 and 4]{Nik11}. We also note that Nikolayevsky's proof for (ii) seems to work in general without nice bases \cite{Nik08ded,Nik11}; thus the novelty of the theorem above is in dropping the assumption on the nice basis in (i) and relating it to the Hilbert--Mumford weight. Theorem \ref{thnik} reduces the problem of whether there exists a nilsoliton to finding a positive solution to the linear equation (\ref{eqnikyyt}), although we first need to determine the finite set $\CF_{\CB} \subset \Fn$ for each orthonormal basis $\CB$ which diagonalises $\phi_{\mu}$, which should be a highly non-trivial problem.

Another special case in which Theorem \ref{thnik} takes a simplified form is when the pre-Einstein derivation has all the eigenvalues simple, since in this case the basis $\CB$ is determined uniquely up to re-ordering. This is exactly the situation to which another theorem of Nikolayevsky \cite{Nik08ded} applies, which we recover as a special case of Theorem \ref{thnik}.

\begin{corollary} \emph{(Nikolayevsky \cite[Theorem 1 and Corollary 1]{Nik08ded})} \label{crnikse}
	Under the assumptions of Theorem \ref{thnik}, suppose furthermore that all the eigenvalues of the pre-Einstein derivation are simple. Then $\CB = \{ e_i \}_{i=1}^n$ is determined uniquely up to re-ordering $e_1 , \dots , e_n$, and the following conditions are equivalent.
	\begin{enumerate}
		\item A homogeneous nilmanifold associated to $(\Fn, \mu , \langle , \rangle)$ admits a left-invariant Ricci soliton.
		\item  $P_0$ lies in the interior of $\mathrm{Conv} (\CF_{\CB})$.
		\item  The equation (\ref{eqnikyyt}) admits a solution $(\alpha_1 , \dots , \alpha_m ) \in \rl^m$ satisfying $\alpha_1 , \dots , \alpha_m >0$ and $\sum_{i=1}^m \alpha_i = \Vert P_0 \Vert^{-2}$.
	\end{enumerate}
\end{corollary}

An important point is that the above corollary covers infinitely many examples, as explained in \cite[\S 3]{Nik08ded}. More specifically, it determines the existence of nilsolitons on certain classes of filiforms, i.e.~nilpotent Lie algebras whose descending central series has the maximal possible length \cite[Theorem 2]{Nik08ded}. See also \cite[Proof of Theorem 9]{Pay10} for another example to which Corollary \ref{crnikse} applies.

As a further application of the variational formalism for the proof of Theorem \ref{thmain}, we propose a significantly modified version of the conjecture by Taketomi--Tamaru \cite{TakTam18} and prove it, providing an obstruction to the existence of nilsolitons; see \S \ref{scmttc} for the comparison to the original conjecture.

\begin{theorem} \emph{(see Theorem \ref{ppttcjw} for details)}
	Let $(N, \textsl{g}_N)$ be a homogeneous nilmanifold and $(\Fn , \mu )$ be the corresponding nilpotent Lie algebra with the inner product $\langle , \rangle$, with respect to which the pre-Einstein derivation is self-adjoint. Let $G$ be a Lie group associated to the endomorphisms of $\Fn$ which commute with $\CS$ and are orthogonal to $\CS$, such that its Lie algebra is written as $\Fp \oplus \Fk$ in terms of the Cartan decomposition, and let $Y$ be the associated symmetric space (see Definitions \ref{dflagk} and \ref{dfgpgk} for the precise statements). Suppose that the following hold:
	\begin{enumerate}
		\item the action of $\mathrm{Aut} (\mu ) \cap G$ on $Y$ is not transitive,
		\item all $\mathrm{Aut} ( \mu ) \cap G$-orbits are congruent to each other with respect to $G$, 
		\item for any $g \in G \setminus  \mathrm{Aut} (\mu ) $, there exists a non-zero element in $\Fp \cap ( \mathrm{Der} (\mu) \setminus \mathrm{Der} (\rho (g) \cdot \mu))$.
	\end{enumerate}
	Then $N$ cannot admit left-invariant Ricci solitons.
\end{theorem}

It would be very interesting to find \textit{new} examples of nilpotent Lie groups that admit no nilsolitons by applying Theorem \ref{thmain}; finding nilpotent Lie groups with no nilsolitons is a non-trivial task, partly because candidates must be sought in dimensions at least 7 (see e.g.~\cite{FerCul14}) since all nilpotent Lie groups of dimension less than or equal to 6 are known to have nilsolitons \cite{Lau02,Wil03}. We only point out that there are infinitely many known examples to which our main results (Theorems \ref{thmain} and \ref{thnik}) apply, via Corollary \ref{crnikse} and \cite[Theorem 2]{Nik08ded} worked out by Nikolayevsky.

\subsection{Related topics} \label{scrtcs}

The relationship between the nilsolitons and GIT stability (or moment maps) was realised by Lauret \cite{Lau01,Lau01b,Lau03,Lau06}. The usual formalism of GIT and moment maps, which is for complex projective varieties, does not naively apply to the situation for nilsolitons, but there are various works which establish similar results for real Lie groups; see e.g.~\cite{BilWin23,BilWin23p,EbeJab09,Mar01,RicSlo90}. The reader is also referred to surveys \cite{Lau09,Jab23} for more details, and \cite{Lau02,Lau10,Lau11,LauWil11,Nik11} for many interesting results on nilsolitons that (at least partially) involve GIT.

The Hilbert--Mumford criterion, which is a fundamental result in GIT \cite{MFK94}, seems to be first employed by Jablonski \cite[Section 6]{Jab11} to find a criterion for the existence of nilsolitons. His result, which the author was unfortunately not aware of when the first version of this paper was posted on arXiv, is similar to Theorem \ref{thmain} but there are two key differences. Firstly, Theorem \ref{thmain} does not need checking the reductivity of the Lie algebra $\Fh_{\phi}$ defined in \cite[Section 6, Step 2]{Jab11}; in the notation of this paper, $\Fh_{\phi}$ equals $\Fg \cap \mathfrak{sl} (\Fn ) \cap \mathrm{Der} (\mu)$. Secondly, the stability condition in Theorem \ref{thmain} is more stringent and hence can be more effective as an obstruction to the existence of nilsolitons, since in \cite[Section 6, Step $3'$]{Jab11} we check $\nu(\lambda  ; \mu )  \ge 0$ for a smaller range of $\lambda$ (orthogonal complement of a subspace strictly larger than $\CS$). We note in passing that \cite[Section 6, Step 3]{Jab11} deals with the case where the inner product is not compatible with $\Fh_{\phi}$ (unlike the setting of this paper), and \cite[Section 6, Step 1]{Jab11} is essentially a review of the pre-Einstein derivation. Finally, we point out that the proof in this paper is different from the one in \cite{Jab11} which relies on the closedness of the group orbit; this paper uses the variational principle and the asymptotic slope of the Kempf--Ness energy, which the author believes is direct and geometric.

There are of course many other related works, some of which are cited above. Nikolayevsky's result \cite[Theorem 2]{Nik11} and his variational formulation have many similarities to the approach in this paper (see e.g.~Remarks \ref{rmdfggphpe} and \ref{rmnkcpknpe}), but a key difference is that we consider the numerical invariant $\nu (\lambda ; \mu )$ while \cite[Theorem 2]{Nik11} is concerned with the closedness of the orbit which can be difficult to check in concrete examples. In the usual situation of GIT for complex projective varieties, the equivalence between the Hilbert--Mumford criterion and the closedness of the orbit is a well-known classical result which holds in general, but the situation seems more subtle for real Lie groups; see e.g.~\cite{BilWin23,BilWin23p,RicSlo90}.

Strictly speaking, the GIT stability that we consider in this paper is the \textit{relative} stability, in which we consider stability that is ``orthogonal'' to directions that are deemed unstable; see Remark \ref{rmrsfrs} for why considering relative stability is necessary in this paper. This notion was introduced by Sz\'ekelyhidi \cite{Szethesis,Sze07}, and played a very important role in the study of extremal K\"ahler metrics. In this paper, we need to consider the stability that is orthogonal to $\CS$, where $\CS$ plays the role of the extremal vector field or the K\"ahler--Ricci soliton vector field. It is important that $\CS$ can be defined irrespectively of the existence of nilsolitons (proved by Nikolayevsky \cite{Nik11}, see Theorem \ref{thnkexuqp} and Definition \ref{dfsvscs}), similarly to the extremal vector field \cite[Corollary D]{FutMab95} or the K\"ahler--Ricci soliton vector field \cite[Lemma 2.2]{TiaZhu02}.

\medskip

\noindent \textbf{Organisation of the paper.} We review preliminary materials in \S \ref{scpmons} which contains no new results. We set up the variational formalism for the proof of Theorem \ref{thmain} in \S \ref{sckne}; it contains a review of known materials but there are also some new definitions and results, particularly in \S \ref{sckneg}, including minor modifications of ideas that appeared in the literature in the past. The final section \S \ref{scpftr} contains the proof of the main results.

\medskip

\noindent \textbf{Acknowledgements} The author thanks Hisashi Kasuya for many helpful discussions, Hiroshi Tamaru for many helpful comments and answering questions, and the anonymous referees for helpful comments. This work is partially supported by JSPS KAKENHI Grant Number JP23K03120 and JP24K00524.

\bigskip

\section{Preliminaries on nilsolitons} \label{scpmons}

\begin{definition}
	A \textbf{homogeneous nilmanifold} $(N, \textsl{g}_N)$ is a pair of a nilpotent Lie group $N$ and a left-invariant metric $\textsl{g}_N$ on $N$.
\end{definition}

A more standard way of defining a homogeneous nilmanifold may be as a connected Riemannian manifold with a transitive and isometric action of a nilpotent Lie group, but a classical result due to Wilson \cite{Wil82} states that such a manifold is isometric to a nilpotent Lie group with a left invariant metric. We assume throughout that $N$ is connected and simply connected.

We recall that an $n$-dimensional homogeneous nilmanifold corresponds one-to-one with a metric nilpotent Lie algebra of dimension $n$.

\begin{definition}
	An $n$-dimensional \textbf{metric Lie algebra} is a triple $(\Fn, \mu, \langle , \rangle )$ consisting of an $n$-dimensional vector space $\Fn$ over $\rl$, the Lie bracket $\mu \in \wedge^2 \Fn^* \otimes \Fn$ satisfying the Jacobi identity, and an inner product $\langle , \rangle$ on $\Fn$.
\end{definition}

Note that the set of Lie brackets on $\Fn$ is a real algebraic subset in $\wedge^2 \Fn^* \otimes \Fn$ cut out by the Jacobi identity. In this paper, which focuses on homogeneous nilmanifolds, we shall always assume that $(\Fn, \mu)$ is nilpotent.

We fix some conventions for the metric Lie algebras that we use throughout this paper. When we choose a basis for $\Fn$ in the rest of the paper, it is always meant to be a $\langle , \rangle$-orthonormal basis unless otherwise stated. $\mathrm{End} ( \Fn )$ stands for the set consisting of endomorphisms on $\Fn$ (over $\rl$), and $GL (\Fn )$ stands for the subset consisting of invertible ones. We consistently write $e^x$ for the scalar exponential function, and $\exp (A)$ for the matrix exponential function for $A \in \mathrm{End} ( \Fn )$. 

We write
\begin{equation*}
	\rho : GL( \Fn ) \to GL ( \wedge^2 \Fn^{*} \otimes \Fn )
\end{equation*}
for the map induced by the natural action; more explicitly, it defines a left action
\begin{equation*}
	( \rho (g) \cdot \mu) (X_1 , X_2) := g \mu (g^{-1} X_1 , g^{-1} X_2 ) ,
\end{equation*}
for any $X_1 , X_2 \in \Fn$ and $g \in GL (\Fn)$. In what follows, we shall often write
\begin{equation*}
	\mu_g := \rho (g) \cdot \mu.
\end{equation*}

\begin{definition}
	A \textbf{derivation} of $\mu$ is an endomorphism $D \in \mathrm{End} ( \Fn )$ which satisfies
\begin{equation*}
	D(\mu (X_1,X_2)) = \mu( D(X_1), X_2 ) + \mu ( X_1, D(X_2))
\end{equation*}
for all $X_1 , X_2 \in \Fn$. The set of all derivations, which is naturally a Lie algebra with respect to the commutator, is written as $\mathrm{Der} (\mu )$. We also write $\mathrm{Aut} (\mu)$ for the connected Lie subgroup of $GL (\Fn)$ corresponding to $\mathrm{Der} (\mu)$, which is the automorphism group of $(\Fn , \mu )$.
\end{definition}

\begin{lemma} \label{lmdvaut}
	For $\lambda \in \mathrm{End} (\Fn)$, $\lambda \in\mathrm{Der} (\mu)$ holds if and only if $\rho( \exp (\lambda t) ) \cdot \mu = \mu$ for all $t \in \rl$.
\end{lemma}

We recall the definition of Ricci solitons, which can be defined for a general Riemannian manifold $(M , \textsl{g})$. We write $\mathrm{Ric}_{\textsl{g}}$ for the Ricci endomorphism of $\textsl{g}$, and $\mathrm{ric}_{\textsl{g}}$ for the Ricci curvature tensor; the relationship between them is such that
\begin{equation} \label{eqrlrcre}
	\textsl{g} ( \mathrm{Ric}_{\textsl{g}}X_1,X_2) = \mathrm{ric}_{\textsl{g}}(X_1,X_2)
\end{equation}
holds for any vector fields $X_1,X_2$ on $M$.

\begin{definition}
	A \textbf{Ricci soliton} is a Riemannian metric $\textsl{g}$ on $M$ satisfying $\mathrm{ric}_{\textsl{g}} = c \textsl{g} + L_X \textsl{g}$, where $L_X$ is the Lie derivative along a vector field $X$ on $M$. We call $X$ above the \textbf{soliton vector field}.
\end{definition}

The Ricci soliton is a self-similar solution to the Ricci flow and plays a very important role in Riemannian geometry. The special case when the soliton vector field is zero is the Einstein metric. In this paper, we study a very special case where the manifold in question is a nilpotent Lie group.

The Ricci curvature of a homogeneous nilmanifold $(N , \textsl{g}_N)$ admits a nice formula in terms of the corresponding metric Lie algebra.

\begin{proposition} \emph{(see e.g.~\cite[equation (8)]{Lau01})} \label{ppfmrcnm} 
	The Ricci curvature tensor of the homogeneous nilmanifold $(N , \textsl{g}_N)$, with the corresponding metric Lie algebra $(\Fn, \mu, \langle , \rangle )$, is given by
	\begin{align*}
		\mathrm{ric}_{\mu , \langle , \rangle} (X_1,X_2) &= \langle \mathrm{Ric}_{\mu , \langle , \rangle} X_1,X_2 \rangle \\
		&=- \frac{1}{2} \sum_{i,j=1}^n \langle \mu (X_1,e_i) , e_j \rangle \langle \mu (X_2,e_i) , e_j \rangle + \frac{1}{4} \sum_{i,j=1}^n \langle \mu (e_i,e_j) , X_1 \rangle \langle \mu (e_i,e_j) , X_2 \rangle
	\end{align*}
	for any $X_1,X_2 \in \Fn$, where $\{ e_k \}_{k=1}^n$ is any orthonormal basis for $\Fn$ with respect to $\langle , \rangle$.
\end{proposition}

In terms of the formula in Proposition \ref{ppfmrcnm}, a left-invariant Ricci soliton on a nilpotent Lie group satisfies the following equation.

\begin{proposition}\emph{(see \cite[Proposition 1.1]{Lau01}, \cite[Theorem 4.2]{Lau09}, and \cite[Theorem 1]{Jab14})} \label{pprstnst}
	A left-invariant metric $\textsl{g}_N$ on a nilpotent Lie group $N$, with the corresponding metric Lie algebra  $(\Fn, \mu , \langle , \rangle )$, is a Ricci soliton if and only if it satisfies
	\begin{equation} \label{eqrstnst}
		\mathrm{Ric}_{\mu , \langle , \rangle} = c I_{\Fn} + D
	\end{equation}
	for some $c \in \rl$ and $D \in \mathrm{Der} (\mu)$.
\end{proposition}

\begin{definition}
	A left-invariant Ricci soliton metric on a nilpotent Lie group is called a \textbf{nilsoliton}. By abuse of terminology, we also say that a metric nilpotent Lie algebra $(\Fn, \mu , \langle , \rangle )$ is a nilsoliton if it satisfies (\ref{eqrstnst}).
\end{definition}

\section{Variational formalism} \label{sckne}

\subsection{Pre-Einstein derivation} \label{scsln}

We start by recalling the following foundational definition by Nikolayevsky.

\begin{definition} (Nikolayevsky \cite[Definition 2]{Nik11})
	A derivation $\phi_{\mu}$ of a Lie algebra $(\Fn , \mu )$ is said to be \textbf{pre-Einstein} if it is semisimple, with all the eigenvalues real, and
	\begin{equation*}
		\mathrm{tr} (\phi_{\mu} \psi ) = \mathrm{tr} (\psi )
	\end{equation*}
	for any $\psi \in \mathrm{Der} (\mu)$.
\end{definition}

This terminology is motivated by the fact that the above equation is satisfied for an Einstein derivation when it exists \cite[\S 2]{Nik11}.

\begin{lemma} \label{lmprdgr}
	If $\phi_{\mu}$ is a pre-Einstein derivation of $(\Fn , \mu )$, then $\phi_{\rho (g) \cdot \mu} := g \phi_{\mu} g^{-1}$ is a pre-Einstein derivation of $(\Fn , \rho(g) \cdot \mu )$.
\end{lemma}

\begin{proof}
	We first note $\mathrm{Der} (\rho (g) \cdot \mu) = g \mathrm{Der} (\mu ) g^{-1}$, and hence any element of $\mathrm{Der} (\rho (g) \cdot \mu)$ can be written as $g \psi g^{-1}$ for some $\psi \in \mathrm{Der} (\mu )$. Thus we get
	\begin{equation*}
		\mathrm{tr} ((g \phi_{\mu} g^{-1}) (g \psi g^{-1})) = \mathrm{tr} (\phi_{\mu} \psi) =  \mathrm{tr} (\psi) =  \mathrm{tr} ( g \psi g^{-1})
	\end{equation*}
	as claimed.
\end{proof}

An important theorem presented below was proved by using the Levi--Mal'tsev decomposition of $\mathrm{Der} (\mu )$ by Nikolayevsky.

\begin{theorem} \emph{(Nikolayevsky \cite[Theorem 1]{Nik11})} \label{thnkexuqp}
	Any Lie algebra admits a pre-Einstein derivation, which is unique up to automorphisms. 
\end{theorem}

\begin{definition}
	We say that a metric nilpotent Lie algebra $(\Fn , \mu , \langle , \rangle )$ is \textbf{compatible with the pre-Einstein derivation} $\phi_{\mu}$ if $\phi_{\mu}$ is self-adjoint with respect to $\langle , \rangle$. By abuse of terminology, we also say that $\langle , \rangle$ is compatible with $\phi_{\mu}$, or that $(\Fn , \mu , \langle , \rangle , \phi_{\mu} )$ is compatible.
\end{definition}

In what follows, we fix a pre-Einstein derivation $\phi_{\mu}$ once and for all. If we choose another one, which must be of the form $\varphi \phi_{\mu} \varphi^{-1}$ for some $\varphi \in \mathrm{Aut} (\mu )$, we replace the inner product $\langle , \rangle$ by $\langle \varphi^{-1} \cdot , \varphi^{-1} \cdot \rangle$, which is compatible with $\varphi \phi_{\mu} \varphi^{-1}$. This is consistent with the behaviour of the Ricci endomorphism, which is
\begin{equation*}
	\mathrm{Ric}_{\mu , \langle \varphi^{-1} \cdot , \varphi^{-1} \cdot \rangle} = \varphi \mathrm{Ric}_{\mu , \langle , \rangle} \varphi^{-1}.
\end{equation*}

We observe that an inner product compatible with the pre-Einstein derivation $\phi_{\mu}$ always exists, since $\phi_{\mu}$ is diagonalisable over $\rl$. Henceforth we fix an inner product $\langle , \rangle$ which is compatible with the pre-Einstein derivation, and set up notational conventions that we use in the rest of this paper. Let $(\Fn, \mu , \langle , \rangle , \phi_{\mu})$ be a compatible metric nilpotent Lie algebra. The inner product $\langle , \rangle$ on $\Fn$ naturally induces the ones on $\mathrm{End}(\Fn)$ and $\wedge^2 \Fn^* \otimes \Fn$. By abuse of notation we denote all these inner products by $\langle , \rangle$. We also note that the inner product on $\mathrm{End}(\Fn)$ induced by $\langle , \rangle$ is the Hilbert--Schmidt inner product, defined for $A, B \in \mathrm{End}(\Fn)$ as $\mathrm{tr} (A^t B)$, where $A^t$ is the $\langle , \rangle$-adjoint of $A$. The fixed inner product $\langle , \rangle$ means that we have a metric isomorphism $\Fn \isom \Fn^*= \mathrm{Hom}_{\rl} (\Fn , \rl )$; in particular, we shall often identify $\mathrm{Sym}^2 (\Fn^*)$ with the set of all self-adjoint endomorphisms on $\Fn$ with respect to $\langle , \rangle$, noting that real symmetric matrices are self-adjoint with respect to the standard inner product. $O(n)$ stands for the isometry group of $\langle , \rangle$. Since we work with the fixed inner product, we shall write $\mathrm{Ric}_{\rho (g) \cdot \mu}$ to denote $\mathrm{Ric}_{\rho (g) \cdot \mu , \langle , \rangle}$. Note finally that $\mathrm{End} (\Fn)$ is naturally a Lie algebra with respect to the commutator $[,]$.

\begin{theorem} \emph{(Nikolayevsky \cite[Theorem 1]{Nik11})} \label{thnikrecpe}
	Suppose that $(\Fn , \mu , \langle , \rangle)$ is compatible with the pre-Einstein derivation $\phi_{\mu}$. If $(\Fn , \mu , \langle , \rangle)$ is a nilsoliton, then its Ricci endomorphism $\mathrm{Ric}_{\mu }$ is proportional to $I_{\Fn} - \phi_{\mu}$.
\end{theorem}

Nikolayevsky \cite[Theorem 1]{Nik11} furthermore proved that all eigenvalues of $\phi_{\mu}$ are all strictly positive if $(\Fn , \mu , \langle , \rangle)$ is a nilsoliton, which is known to be helpful in dealing with concrete examples (see e.g.~\cite{FerCul14}).

Theorem \ref{thnikrecpe} was proved as follows. If $\mathrm{Ric}_{\mu}$ is the Ricci endomorphism of a nilsoliton, then it must be orthogonal to any $\psi \in \mathrm{Der} (\mu )$ with respect to $\langle , \rangle$. This observation was made by Lauret \cite[equation (2)]{Lau01b} for the case when $\psi$ is self-adjoint with respect to $\langle , \rangle$, but inspection of his argument easily implies that it holds for any $ \psi \in \mathrm{Der} ( \mu )$, by computation that is well-known to the experts; see \cite[Proof of Lemma 3.3]{Heb98} and also the proof of Lemma \ref{lmcpnst} presented later. We thus find $\langle \mathrm{Ric}_{\mu } , \psi \rangle = 0$ for any $\psi \in \mathrm{Der} (\mu )$. On the other hand, the same statement holds for $I_{\Fn} - \phi_{\mu}$ since we have
\begin{equation*}
	\langle I_{\Fn} - \phi_{\mu} , \psi \rangle = \mathrm{tr} ((I_{\Fn} - \phi_{\mu})^t \psi ) = \mathrm{tr}( (I_{\Fn} - \phi_{\mu})\psi ) = 0
\end{equation*}
where $\phi_{\mu}^t = \phi_{\mu}$ as $\langle , \rangle$ is compatible with $\phi_{\mu}$. Noting that $\mathrm{Ric}_{\mu }$ and $I_{\Fn} - \phi_{\mu}$ are both contained in $\rl I_{\Fn} \oplus \mathrm{Der} (\mu )$, we find that both of them must be contained in the orthogonal complement of $\mathrm{Der} (\mu )$ inside $\rl I_{\Fn} \oplus \mathrm{Der} (\mu )$ with respect to $\langle , \rangle$, which is of dimension at most 1. Hence $\mathrm{Ric}_{\mu }$ and $I_{\Fn} - \phi_{\mu}$ must be proportional to each other; neither of them are zero since $ I_{\Fn} \not\in \mathrm{Der} (\mu )$ and $\mathrm{Ric}_{\mu}$ cannot be zero by \cite[Theorem 2.4]{Mil76}.

It turns out that the arguments in this paper can be slightly streamlined by considering the following vector space, as opposed to dealing with the fixed pre-Einstein derivation.

\begin{definition} \label{dfsvscs}
	Let $\phi_{\mu}$ be the pre-Einstein derivation of $(\Fn , \mu ) $. We define a vector subspace $\CS$ of $\mathrm{End} (\Fn )$ by
	\begin{equation*}
		\CS := \rl (I_{\Fn} - \phi_{\mu} ). 
	\end{equation*}
\end{definition}

The following result seems to be well-known to the experts, since similar results appear in \cite[\S 6.3 and \S 6.4]{Heb98} and \cite[\S 3]{Nik11} (see also \cite[\S 4]{Jab10}).

\begin{lemma} \label{lmnscmp}
	Suppose that $\langle , \rangle$ is compatible with the pre-Einstein derivation and that there exists $g \in GL(\Fn)$ such that $(\Fn , \rho(g) \cdot \mu , \langle , \rangle )$ is a nilsoliton. Then there exists $h \in GL ( \Fn)$ which commutes with $\phi_{\mu}$ and is self-adjoint and positive definite with respect to $\langle , \rangle$, such that $(\Fn , \rho (h) \cdot \mu , \langle , \rangle )$ is a nilsoliton.
\end{lemma}

\begin{proof}	
	If $(\Fn , \rho(g) \cdot \mu , \langle , \rangle )$ is a nilsoliton, we have
	\begin{equation*}
		\mathrm{Ric}_{\rho (g) \cdot \mu } = c ( I_{\Fn} - g \phi_{\mu} g^{-1} )
	\end{equation*}
	for some $c \in \rl$, noting that $g \phi_{\mu} g^{-1}$ is the pre-Einstein derivation for $(\Fn , \rho (g) \cdot \mu )$ by Lemma \ref{lmprdgr} and Theorem \ref{thnkexuqp} (we replace $g$ by $g \varphi$ for some $\varphi \in \mathrm{Aut}( \mu )$ if necessary). Since the Ricci curvature tensor is symmetric, $\mathrm{Ric}_{\rho (g) \cdot \mu }$ is self-adjoint with respect to $\langle , \rangle$ by (\ref{eqrlrcre}), which implies $(g \phi_{\mu} g^{-1})^t = g \phi_{\mu} g^{-1}$. Since $\phi_{\mu}^t = \phi_{\mu}$, this in turn implies
	\begin{equation*}
		\phi_{\mu} = g^t g \phi_{\mu} (g^t g)^{-1},
	\end{equation*}
	and hence $g^t g$ has to commute with $\phi_{\mu}$. We take $h \in GL (\Fn )$ to be the square root of $g^t g$, so that $h$ is self-adjoint and positive definite with respect to $\langle , \rangle$ with $g = \xi h$ for some $\xi \in O(n)$. We then find $\phi_{\mu} = h^2 \phi_{\mu} h^{-2}$, which implies that $h$ commutes with $\phi_{\mu}$, since $h$ is positive definite and $\phi_{\mu}$ is self-adjoint. Recalling Proposition \ref{ppfmrcnm}, we have
	\begin{equation*}
		\mathrm{Ric}_{\rho (\xi h) \cdot \mu } = \mathrm{Ad}_{\xi } \mathrm{Ric}_{\rho (h) \cdot \mu} 
	\end{equation*}
	and hence
	\begin{equation*}
		\mathrm{Ric}_{\rho (h) \cdot \mu } = c ( I_{\Fn} - h \phi_{\mu} h^{-1} )
	\end{equation*}
	as claimed.
\end{proof}

\begin{lemma} \label{ppsvspps}
	Suppose that the metric nilpotent Lie algebra $(\Fn , \mu , \langle , \rangle )$ is compatible with the pre-Einstein derivation $\phi_{\mu}$. Then we have  $\CS \subset \mathrm{Sym}^2 (\Fn^*)$ and $\dim \CS = 1$. Moreover, for any $g \in GL (\Fn )$ which commutes with $\phi_{\mu}$, the following hold.
	\begin{enumerate}
		\item $\CS = \rl (I_{\Fn} - \phi_{\mu_g} ) \subset ( \rl I_{\Fn} \oplus \mathrm{Der} (\mu_g ) ) \cap \mathrm{Sym}^2 (\Fn^* )$.
		\item $\mathrm{Der} (\mu_g)$ is orthogonal to $\CS$ with respect to $\langle , \rangle$, i.e.~we have $\langle I_{\Fn} - \phi_{\mu} , \psi \rangle = 0$ for any $\psi \in \mathrm{Der} ( \mu_g )$.
	\end{enumerate}
	In particular, for any $g \in GL (\Fn )$ which commutes with $\phi_{\mu}$, we can characterise $\CS$ as the orthogonal complement of $\mathrm{Der} ( \mu_g )$ inside $\rl I_{\Fn} \oplus \mathrm{Der} ( \mu_g )$ with respect to $\langle , \rangle$.
\end{lemma}

\begin{proof}
	The claim $\CS \subset \mathrm{Sym}^2 (\Fn^*)$ is obvious since $\phi_{\mu}$ is self-adjoint with respect to $\langle , \rangle$, and $\dim \CS = 1$ follows easily by noting $ I_{\Fn} \not\in \mathrm{Der} (\mu)$. For the second claim, Lemma \ref{lmprdgr} and the assumption imply $\phi_{\mu_g} = g^{-1} \phi_{\mu} g = \phi_{\mu} \in \mathrm{Der} ( \mu_g)$ and $\phi^t_{\mu} = \phi_{\mu}$. Note moreover that for any $\psi \in \mathrm{Der} ( \mu_g)$ there exists $\psi_0 \in \mathrm{Der} (\mu)$ such that $\psi = g \psi_0 g^{-1}$. We find
	\begin{equation*}
		\langle I_{\Fn} - \phi_{\mu} , \psi \rangle = \mathrm{tr} ( ( I_{\Fn} - \phi_{\mu} )^t \psi ) =  \mathrm{tr} ( ( I_{\Fn} - \phi_{\mu} )g \psi_0 g^{-1} ) = \mathrm{tr} ( ( I_{\Fn} - \phi_{\mu} ) \psi_0  ) = 0,
	\end{equation*}
	as required, since $g$ commutes with $\phi_{\mu}$.
\end{proof}

\subsection{Riemannian symmetric space} \label{scrmssp}

In this section we construct a Riemannian symmetric space which ``parametrises'' the endomorphisms that are orthogonal to $\CS$. A motivation for such a construction, and also a geometric interpretation of $\CS$, can be found later in Remark \ref{rmrsfrs}. We continue with our assumption that $(\Fn , \mu , \langle , \rangle , \phi_{\mu})$ is a compatible metric nilpotent Lie algebra.

\begin{definition} \label{dflagk}
	Recalling $\CS := \rl (I_{\Fn} - \phi_{\mu})$, which we regard as a Lie subalgebra of $\mathrm{End} (\Fn)$ with the trivial Lie bracket, we define
	\begin{equation*}
		\tilde{\Fg} := \{ A \in \mathrm{End} (\Fn) \mid [A,X] =0 \textup{ for any } X \in \CS \} ,
	\end{equation*}
	which is the centraliser of $\CS$ in $\mathrm{End} (\Fn)$, and hence a Lie subalgebra of $\mathrm{End} (\Fn)$. Noting $\CS \subset \tilde{\Fg}$, we write
	\begin{equation*}
		\Fg := \CS^{\perp}
	\end{equation*}
	for the orthogonal complement of $\CS$ inside $\tilde{\Fg}$ with respect to $\langle , \rangle$, and define
	\begin{equation*}
		\Fp := \Fg \cap \mathrm{Sym}^2 (\Fn^*).
	\end{equation*}
	Similarly, we define 
	\begin{equation*}
		\mathfrak{k}:= \{ A \in \Fo (n) \mid [A,X] =0 \textup{ for any } X \in \CS \} ,
	\end{equation*}
	which is a Lie subalgebra of $\Fo (n)$.
\end{definition}

\begin{lemma} \label{lmsmpci}
	The vector space $\Fg = \CS^{\perp}$ is a Lie subalgebra of $\mathrm{End} (\Fn )$ with respect to $[,]$, which is isomorphic to the quotient Lie algebra $\tilde{\Fg} / \CS$ and contains $\Fk$. The pair $(\Fg , \Fk)$ is a Riemannian symmetric pair with respect to the Cartan involution $\theta (A) := -A^t$ inherited from $\mathrm{End} (\Fn)$, and $\Fp$ is the $-1$-eigenspace of $\theta$.
\end{lemma}

\begin{proof}
	For the first claim it suffices to show that, for any $A, B \in \Fg$, $[A,B] \in \tilde{\Fg}$ is orthogonal to $\CS$. Indeed we compute, for any $X \in \CS$,
	\begin{equation*}
		 \langle X, [A,B]  \rangle  = \mathrm{tr} (X^t [A,B]) = \mathrm{tr} (XAB-XBA) = \mathrm{tr} ([X,A]B) + \mathrm{tr}(AXB-XBA)=0
	\end{equation*}
	since $X \in \CS \subset \mathrm{Sym}^2 (\Fn^*)$ and $[X,A] =0$. We also note that elements of $\Fk$ are antisymmetric and hence are automatically orthogonal to $\CS$ with respect to $\langle , \rangle$.
	
	For the second claim we prove that $\theta$ is an involutive automorphism of $\Fg$ whose fixed points are precisely the elements of $\Fk$. It suffices to prove that $A \in \Fg$ implies $-A^t \in \Fg$. Indeed, for any $X \in \CS$ and any $A \in \Fg$ we find
	\begin{equation*}
		[X, (-A)^t] = ([X,A])^t  = 0
	\end{equation*}
	and
	\begin{equation*}
		\langle X , (-A)^t \rangle = - \mathrm{tr} (X^t A^t) = -\mathrm{tr} (AX) =-\mathrm{tr} (XA) = - \langle X, A \rangle =0
	\end{equation*}
	by $X \in \CS \subset \mathrm{Sym}^2 (\Fn^*)$.
\end{proof}

\begin{definition} \label{dfgpgk}
	We define $G$ (resp.~$K$) for the connected Lie subgroup of $GL (\Fn)$ associated to $\Fg$ (resp.~$\Fk$). We write $Y$ for the right coset space
	\begin{equation*}
		Y := K \backslash G.
	\end{equation*}
\end{definition}

The right cosets, as opposed to the left cosets, are more natural in the following discussions (see \S \ref{sckneg}). Note that $K$ is compact, and that $G \subset GL( \Fn )$ and $K \subset O(n)$. A well-known result (see e.g.~\cite[Chapter XI, Theorem 8.6]{KN2}) in the theory of Riemannian symmetric spaces gives us the following result.

\begin{lemma} \label{lmgcdcp}
	The right coset space $Y = K \backslash G$ is a Riemannian symmetric space. $Y$ is a complete Riemannian manifold with respect to the bi-invariant metric, and we have a diffeomorphism
	\begin{equation*}
		\exp : \Fp \isom Y
	\end{equation*}
	given by the exponential map.
\end{lemma}

In particular, any right coset $K \cdot g \in Y$ can be represented (uniquely) by $g = \exp (A)$ for some $A \in \Fp$, which is self-adjoint and positive definite with respect to $\langle , \rangle$.

\begin{lemma} \label{lmgccs}
	If $g \in \exp ( \Fp )$, then $\mathrm{Ric}_{\rho (g) \cdot \mu } \in \tilde{\Fg}$.
\end{lemma}

\begin{proof}
	Since $g$ commutes with $\CS$, we find $\phi_{\mu} = \phi_{\rho (g) \cdot \mu}$ by Lemma \ref{lmprdgr}. Recalling that $\phi_{\mu}$ is self-adjoint with respect to $\langle , \rangle$, we find that $\mathrm{Ric}_{\rho (g) \cdot \mu}$ commutes with $\phi_{\mu} = \phi_{\rho (g) \cdot \mu} \in \mathrm{Der} (\rho (g) \cdot \mu )$ by the result due to Heber \cite[Lemma 2.2]{Heb98}. Thus $\mathrm{Ric}_{\rho (g) \cdot \mu }$ commutes with $\CS$.
\end{proof}

\begin{remark} \label{rmdfggphpe}
	The definition of $G$ in Definition \ref{dfgpgk} and $G_{\phi}$ defined by Nikolayevsky \cite{Nik11} are related as $G \cap SL( \Fn ) = G_{\phi}$; see also \cite[\S 6.3 and \S 6.4]{Heb98}. Lemma \ref{lmgccs} is also contained in \cite[Proof of Lemma 1]{Nik08ded}.
\end{remark}

We recall the following classical result concerning the geodesics in $Y$ (see e.g.~\cite[Chapter XI, Theorem 3.2]{KN2} for more details).
\begin{proposition}
	Any geodesic $\gamma (t) $ in $Y$ emanating from $K \cdot g \in Y$ can be written as $\gamma (t) = K \cdot g \exp (\lambda t)$ for some $\lambda \in \Fp$, which can be extended indefinitely for any $t \in \rl$.
\end{proposition}

Geodesic rays, defined on $[0, + \infty )$, play a particularly important role in this paper.

\begin{definition}
	We write $\{ \gamma_{\lambda} (t) \}_{t \ge 0}$ for a geodesic ray in $Y$ defined by $\lambda \in \Fp$ as
	\begin{equation*}
		\gamma_{\lambda} (t) := K \cdot g \exp (\lambda t /2 ),
	\end{equation*}
	where $\gamma_{\lambda} (0) = K \cdot g \in Y$.
\end{definition}

\subsection{Kempf--Ness energy} \label{sckneg}

Our aim in this paper is to find when there exists an inner product $\langle , \rangle_*$ on $\Fn$ so that $(\Fn , \mu , \langle , \rangle_* )$ is a nilsoliton. We first recall that the inner products on $\Fn$ are parametrised by the right coset space $O(n) \backslash GL ( \Fn )$, which in turn parametrises the left-invariant metrics on $N$. Recall also that any element of $O(n) \backslash GL ( \Fn )$ can be represented by a self-adjoint positive definite endomorphism with respect to $\langle , \rangle$ by the polar decomposition. Thus, finding a nilsoliton inner product $\langle , \rangle_*$ is equivalent to finding a positive definite $g \in GL (\Fn )$ such that $(\Fn, \mu , \langle g \cdot , g \cdot \rangle )$ is a nilsoliton. As explained in \cite[\S 4]{Lau09}, it is isometric to $(\Fn, \rho (g) \cdot \mu , \langle ,\rangle )$ by
\begin{equation*}
	g : (\Fn, \mu , \langle g \cdot , g \cdot \rangle ) \isom (\Fn, \rho (g) \cdot \mu , \langle ,\rangle )
\end{equation*}
and it is more convenient for us to work with the varying Lie bracket and the fixed inner product $\langle , \rangle$ which we assume is compatible with the pre-Einstein derivation, as we can see in the following definition which plays a fundamentally important role in the proof of Theorem \ref{thmain}.

\begin{definition}
	Let $\Vert \cdot \Vert$ be the norm on $\wedge^2 \Fn^* \otimes \Fn$ naturally induced from $\langle , \rangle$ on $\Fn$. We define the \textbf{Kempf--Ness energy} $E_{\mu} : Y \to \rl$ by
	\begin{equation*}
		E_{\mu} (K \cdot g) :=  \log \Vert \rho (g) \cdot \mu \Vert^2 , 
	\end{equation*}
	where we note that $E_{\mu}$ is well-defined since the left $\rho( K )$-action preserves the norm.
\end{definition}

In Lemmas \ref{lmkncvx} and \ref{lmcpnst} presented below, we prove that $(\Fn, \mu , \langle ,\rangle )$ admits a nilsoliton if and only if $E_{\mu}$ has a critical point, and that $E_{\mu}$ is geodesically convex, i.e.~convex along geodesics in $Y$.

\begin{remark} \label{rmknera}
Recalling that the right action of $G $ on $Y$ is an isometry with respect to the bi-invariant metric, we find that the group $\mathrm{Aut} (\mu) \cap G$ acts isometrically on $Y$ from the right. Lemma \ref{lmdvaut} shows that $E_{\mu}$ is invariant under the right $\mathrm{Aut} (\mu) \cap G$-action, and the map $E_{\mu} : Y \to \rl$ naturally descends to the quotient topological space $Y / (\mathrm{Aut} (\mu) \cap G)$.
\end{remark}

\begin{remark} \label{rmrsfrs}
	It is natural to consider the Kempf--Ness energy defined on the right coset space $O(n) \backslash GL(\Fn )$, but the difficulty in this approach is that the energy $ \log \Vert \rho (g) \cdot \mu \Vert^2$ can never be bounded from below over $O(n) \backslash GL(\Fn )$ since multiplying $\mu$ by an overall constant can scale it down to 0 (see also \cite[\S 7]{Lau09}). Such a problem can often be diverted by restricting the domain to $SO(n) \backslash SL ( \Fn )$, which is indeed the case for the original setting of the Kempf--Ness theorem (i.e.~complex projective varieties). This approach fails, however, for nilsolitons, since the space of derivations $\mathrm{Der} (\mu )$ (along which the Kempf--Ness energy should be invariant) is in general not contained in $\mathfrak{sl} (\Fn )$ and it is not translation invariant; if $D$ is a derivation, $c I_{\Fn} + D$ can never be a derivation if $c \neq 0$. Thus the key issue is to single out one specific direction in $\rl I_{\Fn} \oplus \mathrm{Der} (\mu)$, which is orthogonal to $\mathrm{Der} (\mu )$ and along which we cannot hope the Kempf--Ness energy to be bounded from below. We claim that $\CS$ is such a direction, and that we can circumvent the unboundedness problem above by considering the subspace that is orthogonal to $\CS$. Indeed, Lemma \ref{ppsvspps} states that $\CS$ is the orthogonal complement of $\mathrm{Der} (\mu )$ in $\rl I_{\Fn} \oplus \mathrm{Der} (\mu )$ with respect to $\langle , \rangle$, and the Kempf--Ness energy is always unbounded below over $\CS$, since for $X := ( I_{\Fn} - \phi_{\mu})/2 \in \CS$ we have
	\begin{equation*}
		\lim_{t \to + \infty} \frac{\log \Vert \rho ( \exp ( Xt )) \cdot \mu \Vert^2}{t} = \lim_{t \to + \infty} \frac{\log \Vert e^{-t/2} \mu \Vert^2}{t} = -1 <0
	\end{equation*}
	by Lemma \ref{lmdvaut}. The observation above suggests that the directions orthogonal to $\CS$ in an appropriate Riemannian symmetric space are the relevant ones for our problem, and the preparations made in \S \ref{scsln} and \S \ref{scrmssp} are meant to ensure that this approach works; for example, defining $\Fg$ as the orthogonal complement of $\CS$ in $\tilde{\Fg}$ means that $\CS$ and its key properties are preserved when we vary $\mu$ by $g \in G$ (Lemmas \ref{lmprdgr} and \ref{ppsvspps}).
\end{remark}

The lemma below shows that $E_{\mu}$ is geodesically convex, with a specific condition for the strict convexity (see also \cite[Theorem 3.5]{Lau01} for the uniqueness). This result is the key ingredient in the proof of Theorem \ref{thmain}.

\begin{lemma} \label{lmkncvx}
	$E_{\mu}$ is convex along geodesics in $Y$. Moreover, if $E_{\mu}$ is bounded from below, it is strictly convex along geodesics that are not contained in $\mathrm{Aut} (\mu)$ and constant along geodesics contained in $\mathrm{Aut} (\mu)$. If it fails to be strictly convex along a geodesic not contained in $\mathrm{Aut} (\mu)$, then there exists $\lambda \in \Fp$ and a geodesic ray $\{ \gamma_{\lambda} (t) \}_{t \ge 0} \subset Y$ such that
	\begin{equation*}
		\lim_{t \to + \infty} \frac{E_{\mu} (\gamma_{\lambda} (t))}{t} <0.
	\end{equation*}
	
	In particular, if $E_{\mu}$ admits a critical point over $Y$, it must attain the global minimum which is unique up to the $\mathrm{Aut} (\mu) \cap G$-action.
\end{lemma}

\begin{proof}
	Pick $K \cdot g \in Y$ with $g \in \exp ( \Fp )$ as before, which is self-adjoint and positive definite with respect to $\langle , \rangle$. Consider a geodesic $\{ \gamma_{\lambda} (t) \}_{t \ge 0}$ emanating from $K \cdot g$, with $\lambda \in \Fp$. Computing
	\begin{equation*}
		\frac{d}{dt} E_{\mu} (\gamma_{\lambda} (t)) = \frac{1}{\Vert \rho (g \exp (\lambda t/ 2) ) \cdot \mu \Vert^2} \left\langle \rho (g) \rho_* (\lambda) \rho (\exp (\lambda t/ 2) ) \cdot \mu , \rho (g) \rho (\exp (\lambda t/ 2) ) \cdot \mu \right\rangle,
	\end{equation*}
	where
	\begin{equation*}
		\rho_* : \mathrm{End} (\Fn) \to \mathrm{End} (\wedge^2 \Fn^{*} \otimes \Fn)
	\end{equation*}
	is the Lie algebra homomorphism associated to $\rho$, we differentiate the above once more to find
	\begin{align*}
		\left. \frac{d^2}{dt^2} \right|_{t=0} E_{\mu} (\gamma_{\lambda} (t)) = &\frac{1}{2 \Vert \rho (g) \cdot \mu \Vert^2} \left\langle \rho (g) \rho_* (\lambda)^2 \cdot \mu , \rho (g) \cdot \mu \right\rangle + \frac{\Vert \rho (g) \rho_* (\lambda) \cdot \mu \Vert^2}{2 \Vert \rho (g) \cdot \mu \Vert^2}  \\
		&-\frac{1}{\Vert \rho (g) \cdot \mu \Vert^4} \left\langle \rho (g) \rho_* (\lambda) \cdot \mu , \rho (g) \cdot \mu \right\rangle^2 .
	\end{align*}
	Consider an inner product on $\wedge^2 \Fn^* \otimes \Fn$ defined by
	\begin{equation*}
		(a,b)_g:= \left\langle \rho (g) a , \rho (g) b \right\rangle,
	\end{equation*}
	which is indeed well-defined since $g$ is positive definite. We pick an orthonormal basis $\{ \tilde{e}_{\alpha} \}_{\alpha}$ for $\wedge^2 \Fn^* \otimes \Fn$ with respect to this inner product which also diagonalises $\rho_* ( \lambda)$, which is possible since $ \rho_* ( \lambda )$ is self-adjoint with respect to $\langle , \rangle$ (as $\lambda$ is) and $(,)_g$ is positive definite. Presenting $\mu$ with respect to this basis as well, we find
	\begin{align*}
		&\left( \frac{1}{2} \left\langle \rho (g) \rho_* (\lambda)^2 \cdot \mu , \rho (g) \cdot \mu \right\rangle + \frac{1}{2} \Vert \rho (g) \rho_* (\lambda) \cdot \mu \Vert^2 \right) \Vert \rho (g) \cdot \mu \Vert^2 - \left\langle \rho (g) \rho_* (\lambda) \cdot \mu , \rho (g) \cdot \mu \right\rangle^2 \\
		&= \left( \sum_{\alpha} \rho_* ( \lambda)^2_{\alpha \alpha} (\mu_{\alpha})^2  \right) \left( \sum_{\alpha}  (\mu_{\alpha})^2  \right) - \left( \sum_{\alpha} \rho_* ( \lambda)_{\alpha \alpha} (\mu_{\alpha})^2 \right)^2 \ge 0
	\end{align*}
	by the Cauchy--Schwarz inequality for semi-inner products, with equality if and only if there exists $c \in \rl$ such that
	\begin{equation*}
		\sum_{\alpha} (\rho_* ( \lambda)_{\alpha \alpha} -c)^2 (\mu_{\alpha})^2 = 0. 
	\end{equation*}
	The above equation holds if and only if $\rho_* ( \lambda)_{\alpha \alpha} = c$ for all $\alpha$ with $\mu_{\alpha} \neq 0$. Since we took a diagonalising basis for $\rho_* ( \lambda)$, this happens if and only if $\rho( \exp (\lambda )) \cdot \mu = e^{c} \mu$.
	
	Thus we find that, for the geodesic $\gamma_{\lambda} (t) = K \cdot g \exp (\lambda t/ 2)$,
	\begin{equation*}
		\left. \frac{d^2}{dt^2} \right|_{t=0} E_{\mu} (\gamma_{\lambda} (t)) \ge 0
	\end{equation*}
	with equality if and only if there exists $c \in \rl$ such that $\rho( \exp (\lambda t /2)) \cdot \mu = e^{ct} \mu$. The case $c=0$ happens if and only if $\lambda \in \mathrm{Der} (\mu)$ by Lemma \ref{lmdvaut}; note that in this case $E_{\mu} (\gamma_{\lambda} (t))$ is constant in $t$. The case $c \neq 0$ happens only when $E_{\mu} (\gamma_{\lambda} (t))$ is unbounded below with
	\begin{equation*}
		\lim_{t \to + \infty} \frac{E_{\mu} (\gamma_{\lambda} (t))}{t} < 0,
	\end{equation*}
	which is easy to see when $c<0$ (when $c>0$ we replace $\lambda$ by $- \lambda$). In particular, if $E_{\mu}$ is bounded from below, it is strictly convex along geodesics that are not contained in $\mathrm{Aut} (\mu)$, i.e.
	\begin{equation*}
		\left. \frac{d^2}{dt^2} \right|_{t=0} E_{\mu} (\gamma_{\lambda} (t)) >0
	\end{equation*}
	holds for any geodesic $\gamma_{\lambda} (t) = K \cdot  g \exp (\lambda t/ 2)$ with $\lambda \in \Fp \setminus \mathrm{Der} (\mu)$. Moreover, $E_{\mu}$ is constant along any geodesic $\gamma_{\lambda'} (t) = K \cdot  g \exp (\lambda' t/ 2)$ with $\lambda' \in \mathrm{Der} (\mu)$ since $\rho( \exp (\lambda' t /2)) \cdot \mu = \mu$ by Lemma \ref{lmdvaut}.
	
	The convexity of $E_{\mu}$ immediately implies that any critical point must be the global minimum. Thus, if $E_{\mu}$ admits a critical point, it must be bounded from below. In this case $E_{\mu}$ is strictly convex along geodesics not contained in $\mathrm{Aut} (\mu)$, which implies that any two critical points of $E_{\mu}$ must be connected by a geodesic contained in $\mathrm{Aut} (\mu) \cap G$, proving the uniqueness.
\end{proof}

\begin{lemma} \label{lmcpnst}
The Kempf--Ness energy $E_{\mu} : Y \to \rl$ admits a critical point if and only if there exists $g \in GL (\Fn)$ such that $(\Fn , \rho (g) \cdot \mu , \langle , \rangle )$ is a nilsoliton.
\end{lemma}

\begin{proof}
	We first observe that the domain of $E_{\mu}$ can be naturally extended to $G$ to get a map $\tilde{E}_{\mu} : G \to \rl$ defined as
	\begin{equation*}
		\tilde{E}_{\mu} (g ) := \log \Vert \rho (g ) \cdot \mu \Vert^2.
	\end{equation*}
	The first derivative of $\tilde{E}_{\mu}$ along $\lambda \in \Fg$ can be computed as
	\begin{align*}
		\left. \frac{d}{dt} \right|_{t=0} \tilde{E}_{\mu} (g \exp (\lambda t)) &= \frac{2}{\Vert \rho (g) \cdot \mu \Vert^2} \left\langle \rho (g) \rho_* (\lambda) \cdot \mu , \rho (g) \cdot \mu \right\rangle \\
		&=\frac{2}{\Vert \mu_g \Vert^2} \left\langle \rho_* (\mathrm{Ad}_g (\lambda^{\perp} ) ) \cdot \mu_g , \mu_g \right\rangle ,
	\end{align*}
	where we wrote $\lambda^{\perp}$ for the projection of $\lambda$ to the orthogonal complement of $\mathrm{Der} (\mu ) \cap \Fg$ in $\Fg$ with respect to the Hilbert--Schmidt inner product $\langle g^{-1} \cdot , g^{-1} \cdot \rangle$ on $\mathrm{End} (\Fn )$, by noting $\rho_* (D) \cdot \mu =0$ for any $D \in \mathrm{Der} (\mu ) \cap \Fg$ by Lemma \ref{lmdvaut}. Writing $\tilde{\lambda} := \mathrm{Ad}_g (\lambda^{\perp})$, we find
	\begin{align*}
		\left. \frac{d}{dt} \right|_{t=0} \tilde{E}_{\mu} (g \exp (\lambda t)) &= \frac{2}{\Vert \mu_g \Vert^2} \left\langle \rho_* ( \tilde{\lambda} ) \cdot \mu_g , \mu_g \right\rangle \\
		&=\frac{2}{\Vert \mu_g \Vert^2} \left( - \tilde{\lambda}_{\alpha i} (\mu_g )_{\alpha j}^k - \tilde{\lambda}_{\beta j} (\mu_g )_{i \beta}^k + \tilde{\lambda}_{k \gamma} ( \mu_g )_{ij}^{\gamma} \right) (\mu_g )_{ij}^k ,
	\end{align*}
	where we took a $\langle , \rangle$-orthonormal basis $\{ e_i \}_{i=1}^n$ for $\Fn$ and represented $\tilde{\lambda}$ (resp.~$\mu_g$) as a matrix (resp.~a tensor) with respect to this basis, and summation (from $1$ to $n$) is taken over all repeated indices. Exchanging the indices as in \cite[equation (3.3)]{Heb98}, we find
	\begin{align*}
		\left( - \tilde{\lambda}_{\alpha i} (\mu_g )_{\alpha j}^k - \tilde{\lambda}_{\beta j} (\mu_g )_{i \beta}^k + \tilde{\lambda}_{k \gamma} ( \mu_g )_{ij}^{\gamma} \right) (\mu_g )_{ij}^k &= \tilde{\lambda}_{\alpha i} \left( -( \mu_g )_{\alpha j}^k ( \mu_g )_{i j}^k - ( \mu_g )_{j \alpha }^k ( \mu_g )_{ji}^k + ( \mu_g )_{kj}^{\alpha} ( \mu_g )_{kj}^i \right) \\
		&= \tilde{\lambda}_{\alpha i} \frac{\mathrm{ric}_{\mu_g} (e_i,e_{\alpha})}{4} \\
		&= \frac{1}{4} \left\langle \mathrm{Ric}_{\mu_g} , \tilde{\lambda} \right\rangle ,
	\end{align*}
	by Proposition \ref{ppfmrcnm}.
	
	Suppose now that $E_{\mu}$ admits a critical point $K \cdot g \in Y$. By Lemma \ref{lmkncvx}, any critical point of $E_{\mu}$ must be the global minimum, which in turn has to minimise $\tilde{E}_{\mu}$ over $G$. Thus, $g \in G$ needs to attain the global minimum of $\tilde{E}_{\mu}$, and hence it must be the critical point of $\tilde{E}_{\mu}$ satisfying
	\begin{equation*} 
		\left\langle \mathrm{Ric}_{\mu_g} , \mathrm{Ad}_g (\lambda^{\perp}) \right\rangle = 0 
	\end{equation*}
	for any $\lambda \in \Fg$. This condition holds if and only if we have
	\begin{equation*}
		\mathrm{Ric}_{\mu_g} \in \CS \oplus g \mathrm{Der} (\mu ) g^{-1} = \CS \oplus \mathrm{Der} (\mu_g) \subset \rl I_{\Fn} \oplus \mathrm{Der} (\mu_g) ,
	\end{equation*}
	since $\mathrm{Ric}_{\mu_g} \in \tilde{\Fg}$ by Lemma \ref{lmgccs} and $\mathrm{Ad}_g (\lambda^{\perp}) \in \Fg$ is orthogonal to $\mathrm{Der} (\mu_g)$ with respect to $\langle , \rangle$, where we recall that $\lambda^{\perp}$ is orthogonal to $\mathrm{Der} (\mu)$ with respect to the Hilbert--Schmidt inner product  $\langle g^{-1} \cdot , g^{-1} \cdot \rangle$ on $\mathrm{End} (\Fn )$, i.e.~$\mathrm{tr} ((g \lambda^{\perp} g^{-1})^t g D g^{-1})=0$ for all $D \in \mathrm{Der} (\mu)$. Recall also that $\CS \oplus \mathrm{Der} (\mu_g) \subset \rl I_{\Fn} \oplus \mathrm{Der} (\mu_g)$ by Lemma \ref{ppsvspps}. This implies that $(\Fn , \mu_g , \langle , \rangle )$ is a nilsoliton by Proposition \ref{pprstnst}.

	Suppose conversely that there exists $g \in GL (\Fn)$ such that $(\Fn , \rho (g) \cdot \mu , \langle , \rangle )$ is a nilsoliton. By Lemma \ref{lmnscmp}, we may assume without loss of generality that $g$ commutes with $\CS$, and that $g$ is self-adjoint and positive definite with respect to $\langle , \rangle$. This implies that there exists $A \in \tilde{\Fg} \cap \mathrm{Sym}^2 ( \Fn^* ) = \Fp \oplus \CS$ such that $g = \exp ( A )$. Writing $A = B + X$ for some $B \in \Fp$ and $X \in \CS$, we define $h := \exp (B)$ and $\varphi := \exp (X)$ to find $g = h \varphi$ as $[B,X]=0$. We have thus found $h \in \exp ( \Fp )$ such that $(\Fn , \rho (h) \cdot \mu , \langle , \cdot \rangle ) = (\Fn , \rho (h \varphi ) \cdot \mu , \langle , \cdot \rangle )$ is a nilsoliton, where we note $\mathrm{Ric}_{\rho (h) \cdot \mu} \in \CS$ by Theorem \ref{thnikrecpe}. Repeating the computation above for $\lambda \in \Fp$, we find that
	\begin{equation*}
		\left. \frac{d}{dt} \right|_{t=0} E_{\mu} (K \cdot h \exp (\lambda t)) = \frac{1}{2 \Vert \mu_h \Vert^2} \left\langle \mathrm{Ric}_{\rho (h) \cdot \mu} , \mathrm{Ad}_h (\lambda^{\perp}) \right\rangle =0
	\end{equation*}
	for any $\lambda \in \Fp$, since $\mathrm{Ad}_h (\lambda^{\perp}) \in \Fg$ and $\mathrm{Ric}_{\rho (h) \cdot \mu} \in \CS$ are orthogonal to each other by the definition of $\Fg$. Thus $K \cdot h\in Y$ is a critical point of $E_{\mu} : Y \to \rl$ as claimed.
\end{proof}	
	
\begin{remark} \label{rmnkcpknpe}
	Similar results also appear in \cite{Nik11,Heb98}.
\end{remark}

We observe the following consequence of the convexity, which shows that $E_{\mu}$ admits no critical points if $\mu$ is an eigenvector with non-zero eigenvalue of $\rho_* (\lambda)$; the argument is essentially contained in the proof of Lemma \ref{lmkncvx}.

\begin{proposition} \label{ppmevnex}
	Suppose that there exists $\lambda \in \Fp$ such that $\rho ( \exp (\lambda ) ) \cdot \mu = e^{-c} \cdot \mu$ for some $c \neq 0$. Then $E_{\mu}$ admits no critical points over $Y$.
\end{proposition}

\begin{proof}
	If $c>0$, we easily find that $E_{\mu}$ cannot be bounded from below (hence no global minimum exists) since
	\begin{equation*}
		\lim_{t \to + \infty} \frac{E_{\mu} (\gamma_{\lambda} (t))}{t} = -c <0.
	\end{equation*}
	If $c<0$ we replace $\lambda$ by $-\lambda$ and argue as above.
\end{proof}

\section{Proof of the results} \label{scpftr}

We continue with our assumption that $(\Fn , \mu , \langle , \rangle , \phi_{\mu})$ is a compatible metric nilpotent Lie algebra.

\subsection{Hilbert--Mumford criterion}

We prove that the asymptotic slope of the Kempf--Ness energy agrees with the Hilbert--Mumford weight, and the following argument is well-known in Geometric Invariant Theory for complex reductive Lie groups acting on complex projective varieties (see e.g.~\cite{MFK94}). Suppose that we take $\lambda \in \mathrm{Sym}^2 (\Fn^*)$ and diagonalise it as $\mathrm{diag} (\lambda_1 , \dots , \lambda_n)$, by choosing a $\langle , \rangle$-orthonormal basis $\{ e_i \}_{i=1}^n$ for $\Fn$. Suppose that we also write $\mu (e_i, e_j) = \sum_{k=1}^n \mu^k_{ij}e_k$ with respect to this basis.

\begin{definition} \label{dfhmwt}
	The \textbf{Hilbert--Mumford weight} of $\lambda \in \mathrm{Sym}^2 (\Fn^*)$ with respect to the Lie bracket $\mu$ is defined by
	\begin{equation*}
		\nu(\lambda  ; \mu ) := - \min_{1 \le i < j \le n, 1 \le k \le n} \{ \lambda_i+\lambda_j - \lambda_k \mid \mu^k_{ij} \neq 0\}.
	\end{equation*}
\end{definition}

\begin{remark}
While the above definition depends on the $\langle , \rangle$-orthonormal basis for $\Fn$, the Hilbert--Mumford weight can be defined independently of a basis chosen. Indeed, when we choose $\lambda \in \mathrm{Sym}^2 (\Fn^*)$ and decompose $\wedge^2 \Fn \otimes \Fn^* =:V$ in terms of $\rho_* (\lambda )$-eigenspaces as
\begin{equation*}
	V = \bigoplus_{i} V_{ \lambda }^{s_i}
\end{equation*}
where $\{ s_i \}_i \subset \rl$ is the set of eigenvalues of $\rho_* (\lambda )$ with the corresponding eigenspaces $\{ V_{ \lambda }^{s_i} \}_i$, we have a natural $\rl$-filtration of $V$, i.e.~a family of subspaces $\{ F_{\lambda}^{s} V \}_{s \in \rl}$ of $V$ defined by
\begin{equation*}
F^{s}_{\lambda } V := \bigoplus_{s_i \ge s} V_{ \lambda }^{s_i}.
\end{equation*}
It is well-known \cite[\S 1.1]{BHJ1} that a filtration of $V$ naturally corresponds to a non-Archimedean metric on $V$ defined as
\begin{equation*}
	\Vert v \Vert_{\mathrm{NA}( \lambda )} := e^{- \sup \{ s \in \rl \mid v \in F^{s}_{\lambda} V \} }, \quad v \in V.
\end{equation*}
It is straightforward to check by direct computation that
\begin{equation*}
	\nu(\lambda  ; \mu ) = \log \Vert \mu \Vert_{\mathrm{NA}( \lambda )},
\end{equation*}
which gives a basis-independent definition of the Hilbert--Mumford weight. 
\end{remark}

\begin{lemma} \label{lmasemhm}
	Let $\{ \gamma_{\lambda} (t) \}_{t \ge 0}$ be a geodesic ray in $Y$ emanating from $K \cdot I_{\Fn}$, defined by $\lambda \in \Fp$. Then we have
	\begin{equation*}
		\lim_{t \to + \infty} \frac{d}{dt} E_{\mu} (\gamma_{\lambda} (t)) = \lim_{t \to + \infty} \frac{E_{\mu} (\gamma_{\lambda} (t))}{t} = \nu(\lambda  ; \mu ).
	\end{equation*}
\end{lemma}

\begin{proof}
	The first equality is a straightforward consequence of the fact that a difference quotient of a convex function is monotonically increasing \cite[Theorem 1.1.5 and Corollary 1.1.6]{Hor07}. Noting $\lambda \in \Fp \subset \mathrm{Sym}^2 (\Fn^*)$, we compute
	\begin{align*}
		E_{\mu} (\gamma_{\lambda} (t)) &= \log \left( \sum_{k=1}^n \sum_{1 \le i < j \le n}  e^{(\lambda_k  - \lambda_i - \lambda_j )t } (  \mu^k_{ij} )^2 \right)  \\
		&=\nu(\lambda  ; \mu )t + \log \left( \sum_{k=1}^n \sum_{1 \le i < j \le n}  e^{-\nu(\lambda  ; \mu )t-( \lambda_i + \lambda_j - \lambda_k )t } (  \mu^k_{ij} )^2 \right) ,
	\end{align*}
	where it is important to note that the second term above remains bounded as $t \to + \infty$; the leading term inside the $\log$ is $(  \mu^k_{ij} )^2 \neq 0$ whose index attains the minimum in the definition of $\nu (\lambda  ; \mu )$, and other terms decay exponentially. The claimed result follows immediately.
\end{proof}

We recall the following theorem.

\begin{theorem} \emph{(see e.g.~\cite[Theorem 4.12]{HasHaya})} \label{thcvhaya}
	Let $(Z,\textsl{g}_Z)$ be a complete Riemannian manifold with a fixed base point $z_0 \in Z$, and $H$ be a Lie group acting isometrically on $(Z,\textsl{g}_Z)$ from the right. Let $f : Z \to \rl$ be a smooth $H$-invariant function which satisfies
	\begin{equation} \label{eqstcvmh}
		f(\gamma (t) ) \le \frac{t - t_0}{t_1 - t_0}  f(\gamma (t_1)) + \frac{t_1 - t}{t_1 - t_0} f ( \gamma (t_0) )
	\end{equation}
	for any geodesic segment $\{ \gamma (t) \}_{t_0 \le t \le t_1}$ and for any $t \in [t_0 , t_1]$, with equality at some $t \in [t_0 , t_1]$ if and only if $\{ \gamma (t ) \}_{t_0 \le t \le t_1}$ is contained in $\gamma (t_0) \cdot H$. Then the following are equivalent.
	\begin{enumerate}
		\item There exists $z_* \in Z$ which attains the global minimum of $f$, and it is unique modulo the $H$-action.
		\item For any geodesic ray $\{ \gamma (t) \}_{t \ge 0}$ emanating from $z_0$, not contained in the $H$-orbit, we have
		\begin{equation*}
			\lim_{t \to + \infty} \frac{f(\gamma( t ))}{ t } >0.
		\end{equation*}
	\end{enumerate}
\end{theorem}

\begin{remark} \label{rmstcvh}
	The inequality (\ref{eqstcvmh}) is equivalent to saying that $f$ is convex along geodesics, and the condition for the equality states that $f$ is strictly convex along geodesics that are not contained in an $H$-orbit. Indeed, the equality in (\ref{eqstcvmh}) holds for some $t \in [t_0 , t_1]$ if and only if it does so for all $t \in [t_0 , t_1]$ by the convexity, which is equivalent to saying $\frac{d^2 f}{dt^2} (\gamma (t) ) = 0$ for all $t \in [t_0 , t_1]$.
\end{remark}

We now get the following corollary which immediately implies Theorem \ref{thmain}.

\begin{corollary} \label{crkntns}
	The Kempf--Ness energy $E_{\mu} : Y \to \rl$ admits a critical point, which is unique modulo the $\mathrm{Aut} (\mu) \cap G$-action, if and only if
	\begin{equation*}
		\lim_{\tau \to + \infty} \frac{E_{\mu}(\gamma_{\lambda} (t))}{t} >0
	\end{equation*}
	holds for any geodesic ray $\{ \gamma_{\lambda} (t) \}_{t \ge 0}$ with $\lambda \in \Fp \setminus \mathrm{Der} (\mu)$.
	
	In particular, a homogeneous nilmanifold $(N, \textsl{g}_N)$ with the metric nilpotent Lie algebra $(\Fn, \mu , \langle , \rangle)$ admits a left-invariant Ricci soliton if and only if 
	\begin{equation*}
	\nu (\lambda  ; \mu ) \ge 0	
	\end{equation*}
	for any $\lambda \in \Fp$, with equality if and only if $\lambda \in \mathrm{Der} (\mu) \cap \Fp$.
\end{corollary}

\begin{proof}
	For the first statement, we apply Theorem \ref{thcvhaya} to $Z = Y$, $\textsl{g}_Z=\text{bi-invariant metric on $Y$}$, $H= \mathrm{Aut} (\mu) \cap G$, and $f = E_{\mu}$; recall that $E_{\mu}$ is invariant under the right $\mathrm{Aut} (\mu) \cap G$-action by Lemma \ref{lmdvaut} (see also Remark \ref{rmknera}). As we saw in the previous section, particularly Lemma \ref{lmkncvx}, the Kempf--Ness energy $E_{\mu}$ satisfies all the hypotheses in Theorem \ref{thcvhaya} when it is strictly convex along geodesics not contained in $\mathrm{Aut} (\mu) \cap G$ (see also Remark \ref{rmstcvh}). When the strict convexity fails for a geodesic not contained in $\mathrm{Aut} (\mu) \cap G$, there exists a geodesic ray $\{ \gamma_{\lambda} (t) \}_{t \ge 0}$ of negative asymptotic slope, again by Lemma \ref{lmkncvx}, and $E_{\mu}$ admits no global minimum as it is unbounded below; thus the claimed equivalence holds in this case as well. The second statement follows from Lemmas \ref{lmkncvx}, \ref{lmcpnst}, and \ref{lmasemhm}, where we also note $\nu (\lambda  ; \mu ) = 0$ for any $\lambda \in \Fp \cap \mathrm{Der} (\mu)$ by the $\mathrm{Aut} (\mu) \cap G$-invariance of $E_{\mu}$.
\end{proof}

The condition for $\nu (\lambda  ; \mu )$ in Corollary \ref{crkntns} is the nilsoliton version of the Hilbert--Mumford criterion. As is well-known in GIT, this criterion can be written equivalently as follows.

\begin{proposition}
	For each $\lambda \in \mathrm{Sym}^2 (\Fn^*)$, define 
	\begin{equation*}
		\CT := \{ (i,j,k) \in \itg^3 \mid 1 \le i < j \le n, 1 \le k \le n, \textup{ s.t.~} \mu^k_{ij} \neq 0 \}
	\end{equation*}
	where we presented $\mu$ with respect to a $\langle , \rangle$-orthonormal basis which diagonalises $\lambda$. Suppose that we re-order the elements of $\{ \lambda_i+\lambda_j - \lambda_k \}_{(i,j,k) \in \CT }$ as $\beta_1 (\lambda) \le \cdots \le \beta_m (\lambda)$, where $m:= | \CT |$. Then the following are equivalent.
	\begin{enumerate}
		\item $\nu (\lambda  ; \mu ) \ge 0$ for any $\lambda \in \Fp$, with equality if and only if $\lambda \in \Fp \cap \mathrm{Der} (\mu)$.
		\item $\beta_1 (\lambda) < 0 < \beta_m (\lambda)$ for any $\lambda \in \Fp \setminus \mathrm{Der} (\mu)$, and $\beta_1 (\lambda) = \cdots = \beta_m (\lambda) = 0$ for any $\lambda \in \Fp \cap \mathrm{Der} (\mu)$.
	\end{enumerate}
\end{proposition}

\begin{proof}
	Obvious by noting $\nu (\lambda  ; \mu ) = - \beta_1 (\lambda)$ and $\nu (- \lambda  ; \mu ) = \beta_m (\lambda)$.
\end{proof}

Note that the above proposition and Corollary \ref{crkntns} immediately recovers Proposition \ref{ppmevnex}. We finally recall that $E_{\mu}$ cannot admit any critical points if $\phi_{\mu} = 0$ by Theorem \ref{thnikrecpe} and \cite[Theorem 2.4]{Mil76}. Corollary \ref{crkntns} implies that in this case there exists $\lambda \in \Fp \setminus \mathrm{Der} ( \mu )$ with $\nu (\lambda ; \mu ) \le 0$, but we can explicitly write down such an endomorphism.

\begin{proposition} \label{lmslmid}
	Suppose that we have $\phi_{\mu} = 0$, or equivalently $\CS = \rl I_{\Fn}$, and that we take $e_n \in \Fn$ of unit length from the last term of the descending central series so that $\mu(e_n, A) =0$ for any $A \in \Fn$, and extend it to a $\langle , \rangle$-orthonormal basis $\{ e_i \}_{i=1}^n$ for $\Fn$. Let $\lambda \in \mathrm{End} (\Fn )$ be presented as $\mathrm{diag} (1, \dots , 1,  -n+1)$ with respect to this basis. Then $\lambda \in \Fp \setminus \mathrm{Der} ( \mu )$ and $\nu (\lambda ; \mu ) < 0$.
\end{proposition}

\begin{proof}
	First note that $\Fg = \mathfrak{sl} (\Fn )$ when $\CS = \rl I_{\Fn}$, and that any inner product $\langle , \rangle$ is compatible with $\phi_{\mu} = 0$. Then we have $\lambda \in \Fp$ as it is symmetric and trace-free. We note
	\begin{equation*}
		\lambda_i + \lambda_j - \lambda_k = \begin{cases}
			1 \quad &( 1 \le i,j,k \le n-1 ) \\
			n+1 \quad &(1 \le i,j \le n-1 ,\; k=n) 
	\end{cases}
	\end{equation*}
	which exhausts all the possibilities for the indices $(i,j,k)$ with potentially non-vanishing $\mu_{ij}^k$, crucially because $\mu_{nj}^k=\mu_{in}^k =0$ for any $1 \le i,j,k \le n$. We get $\nu (\lambda ; \mu ) \le -1 < 0$, which in turn implies $\lambda \not\in \mathrm{Der} ( \mu )$.
\end{proof}

\subsection{Generalisation of Nikolayevsky's criterion} \label{scgncr}

The proof of Theorem \ref{thnik} is based on a well-known argument in GIT (see e.g.~\cite[Theorem 9.2]{Dol}) for the torus action using the theory of convex sets, but a non-trivial modification is necessary since we need to argue relative to $\CS$.

\begin{proof}[Proof of Theorem \ref{thnik}]
	We fix a $\langle , \rangle$-orthonormal basis $\CB = \{ e_i \}_{i=1}^n$ for $\Fn$ which diagonalises $\phi_{\mu}$, so that
	\begin{equation*}
		I_{\Fn} - \phi_{\mu} = \mathrm{diag} (s_1 , \dots , s_n)
	\end{equation*}
	with respect to $\CB$, and write $\CF$ for $\CF_{\CB}$ and $\mathsf{L}$ for $\mathsf{L}_{\CB}$ to simplify the notation. Since $\rho (\exp ( t (I_{\Fn} - \phi_{\mu}) ) ) \cdot \mu = e^{-t} \mu$ by Lemma \ref{lmdvaut}, we find that 
	\begin{equation} \label{eqsijkc}
		- s_i - s_j + s_k = -1
	\end{equation}
	for any $1 \le i,j,k \le n$ with $\mu_{ij}^k \neq 0$.

	Let $\mathsf{S}$ be a line in $\Fn$ generated by $\CS$ as
	\begin{equation*}
		\mathsf{S}:= \left. \left\{\sum_{i=1}^n \tau s_i e_i \; \right| \; \tau \in \rl \right\}.
	\end{equation*}
	We prove that $\mathsf{S}$ is orthogonal to $\mathsf{L}$, and hence $P_0 \in \mathsf{S}$. Indeed, for any $e_i + e_j - e_k \in \CF$, (\ref{eqsijkc}) implies
	\begin{equation} \label{eqdftsot}
		\left\langle e_i + e_j - e_k - \sum_{l=1}^n \tau s_l e_l , \sum_{m=1}^n \tau s_m e_m \right\rangle = \tau (s_i + s_j - s_k) - \tau^2 \sum_{l=1}^n s_l^2 =  \tau \left( 1 - \tau \sum_{l=1}^n s_l^2 \right),
	\end{equation}
	which vanishes for all indices $(i,j,k)$ with $\mu_{ij}^k \neq 0$ when we choose $\tau$ to be equal to
	\begin{equation} \label{eqdftsot2}
		\tau^* :=  \frac{1}{\sum_{l=1}^n s_l^2}  ,
	\end{equation}
	which is possible since $I_{\Fn} - \phi_{\mu} \neq 0$. Note that the above implies $P_0 =  \sum_{i=1}^n s_i e_i / (\sum_{l=1}^n s_l^2 )$ and hence we have
	\begin{equation} \label{eqtst}
		\Vert P_0 \Vert^2 = \frac{1}{\sum_{l=1}^n s_l^2} =  \tau^* \neq 0 .
	\end{equation}

	We now work with the affine hyperplane $\mathsf{S}^{\perp}$ of $\Fn$ which is orthogonal to $\mathsf{S}$ and contains $\mathsf{L}$, where $P_0 = \sum_{i=1}^n \tau^* s_i e_i$ is the origin of this hyperplane. For each $\lambda := \mathrm{diag} (\lambda_1 , \dots , \lambda_n ) \in \rl^n$ with $\sum_{i=1}^n \lambda_i s_i =0$, we define an affine linear map $f_{\lambda} : \Fn \to \rl$ by
	\begin{equation*}
		f_{\lambda} (w) := \langle w - P_0 , \lambda \rangle ,
	\end{equation*}
	which defines a hyperplane $\mathsf{H}_{\lambda} := \{ w \in \Fn \mid f_{\lambda} (w) =0 \} \subset \Fn$ containing $\mathsf{S}$ (and hence contains $P_0$ and $0 \in \Fn$). Thus the subspace
	\begin{equation*}
		\mathsf{H}_{\lambda} \cap \mathsf{S}^{\perp} = \{ w \in \mathsf{S}^{\perp} \mid f_{\lambda} (w) =0 \}.
	\end{equation*}
	is a proper hyperplane in $\mathsf{S}^{\perp}$ passing through $P_0$, since $\mathsf{H}_{\lambda}$ and $\mathsf{S}^{\perp}$ intersect transversally. Note that $\mathsf{H}_{\lambda} \cap \mathsf{L}$ defines a hyperplane in $\mathsf{L}$ if and only if $\mathsf{L} \not\subset \mathsf{H}_{\lambda}$. Conversely, any hyperplane $\mathsf{H}$ in $\mathsf{L}$ passing through $P_0$ can be realised as above; we simply extend $\mathsf{H}$ to an affine hyperplane in $\Fn$ passing through $P_0$ that contains $\mathsf{S}$, which can be written as $\mathsf{H}_{\lambda}$ for some $\lambda = \mathrm{diag} (\lambda_1 , \dots , \lambda_n ) \in \rl^n$ satisfying $\sum_{i=1}^n \lambda_i s_i =0$.
	
	Henceforth we assume $\dim \mathrm{Conv} ( \CF) \ge 1$ and $\Fp_{\CB} \setminus \mathrm{Der} (\mu) \neq \emptyset$, since otherwise the claimed result holds trivially as follows. We first observe that $\Fp_{\CB} \subset \mathrm{Der} (\mu)$ holds if and only if $f_{\lambda} (e_i + e_j - e_k ) = \lambda_i + \lambda_j - \lambda_k = 0$ for any $e_i + e_j - e_k \in \CF$ and any $\lambda = \mathrm{diag} (\lambda_1 , \dots , \lambda_n ) \in \rl^n$ satisfying $\sum_{i=1}^n \lambda_i s_i =0$, which holds if and only if $\mathsf{L} \subset \mathsf{H}_{\lambda}$ for any such $\lambda$. This in turn is equivalent to $\mathsf{L} \subset \mathsf{S}$, which happens if and only if $\mathsf{L}$, and hence $\CF$, is a singleton set $\{ P_0 \}$ by $\mathsf{L} \subset \mathsf{S}^{\perp}$ proved above. We thus conclude that $\Fp_{\CB} \subset \mathrm{Der} (\mu)$ holds if and only if $\dim \mathrm{Conv} ( \CF) =0$, for which we have $\mathsf{L} = \mathrm{Conv} ( \CF) = \CF = \{ P_0 \}$. In this case the claimed result holds by decreeing the interior of $\mathrm{Conv} ( \CF) = \{ P_0 \}$ to be $\{ P_0 \}$ itself.

	Suppose now that $P_0$ is not contained in the interior of $\mathrm{Conv} ( \CF)$. Then there exists an affine linear function $f : \mathsf{L} \to \rl$ such that
	\begin{equation} \label{eqshpth}
	f(P_0)=0 \le f(w)	
	\end{equation}
	for any $w \in \mathrm{Conv} ( \CF)$, and $f(w) >0$ if $w$ is in the interior of $\mathrm{Conv} ( \CF)$ by the supporting hyperplane theorem (see e.g.~\cite[Corollaries 2.1.11 and 2.1.12]{Hor07}). As discussed above, there exists $\lambda = \mathrm{diag} (\lambda_1 , \dots , \lambda_n ) \in \rl^n$ satisfying $\sum_{i=1}^n \lambda_i s_i =0$ such that $f  = f_{\lambda} |_{\mathsf{L}}$, since $f (P_0)=0$. Now, for any $e_i + e_j - e_k \in \CF$, (\ref{eqshpth}) implies
	\begin{equation} \label{eqfllijk}
		f_{\lambda} (e_i + e_j - e_k ) = \langle e_i + e_j - e_k - P_0 , \lambda \rangle = \lambda_i + \lambda_j - \lambda_k \ge 0 .
	\end{equation}
	We prove $\lambda \in \Fp_{\CB} \setminus \mathrm{Der} (\mu )$. Indeed, $\lambda$ is self-adjoint with respect to $\langle , \rangle$ and commutes with $\CS$ as it is diagonal, and $\lambda \in \Fg$ since $\sum_{i=1}^n \lambda_i s_i =0$. Since $f_{\lambda} (w) >0$ if $w$ is in the interior of $\mathrm{Conv} ( \CF)$, there has to exist some $e_i + e_j - e_k \in \CF$ such that $f_{\lambda} (e_i + e_j - e_k )  = \lambda_i + \lambda_j - \lambda_k > 0$, implying $\lambda \not\in \mathrm{Der} (\mu )$. We thus get $\lambda \in \Fp_{\CB} \setminus \mathrm{Der} (\mu )$ which satisfies $\nu (\lambda ; \mu ) \le 0$.

	Suppose conversely that there exists $\lambda \in \Fp_{\CB} \setminus \mathrm{Der} (\mu )$, presented as $\lambda = \mathrm{diag} (\lambda_1 , \dots , \lambda_n )$ with respect to $\CB$, such that $\nu (\lambda ; \mu ) \le 0$. We then find $f_{\lambda} (w) \ge 0$ for any $w \in \mathrm{Conv} ( \CF)$ exactly as in (\ref{eqfllijk}). Note that in this case there exists at least one index $(i,j,k)$ with $\mu_{ij}^k \neq 0$ such that $\lambda_i+ \lambda_j - \lambda_k >0$, since otherwise we have $\lambda \in \mathrm{Der} (\mu)$ which is a contradiction. Thus there exists at least one point $w \in \mathrm{Conv} ( \CF)$ such that $f_{\lambda} (w) >0$, and hence the hyperplane $\mathsf{H}_{\lambda} \cap \mathsf{S}^{\perp}$ cannot contain $\mathsf{L}$ and thus $\mathrm{Conv} ( \CF)$ is contained in the supporting half-space of $f_{\lambda}$ inside $\mathsf{L}$ at $P_0$, implying that $P_0$ cannot be contained in the interior of $\mathrm{Conv} ( \CF)$.
	
	To prove the second statement, we first fix a linear isomorphism $\Fn \isom \rl^n$ by sending $e_i$ to the standard basis vector $\mathbf{e}_i$ for $\rl^n$ for $i=1 , \dots , n$. We observe that $P_0$ is in the interior of $\mathrm{Conv} (\CF )$ if and only if $P_0$ can be written as a strict convex combination of elements in $\CF$. With respect to the given identification $\Fn \isom \rl^n$, we have $P_0 = \sum_{i=1}^n \tau^* s_i \mathbf{e}_i$, which is contained in the interior of $\mathrm{Conv} (\CF )$ if and only if there exist $\beta_1 , \dots , \beta_m >0$ with $\sum_{i=1}^m \beta_i = 1$ such that
	\begin{equation} \label{eqnlae1}
		P_0 = \begin{pmatrix}
			\tau^*  s_1 \\
			\vdots \\
			\tau^*  s_n 
		\end{pmatrix}
		= \mathbf{Y}^t
		\begin{pmatrix}
			\beta_1 \\
			\vdots \\
			\beta_m
		\end{pmatrix}
	\end{equation}
	where we wrote $\mathbf{Y}$ for $\mathbf{Y}_{\CB}$. We also observe that (\ref{eqsijkc}) implies
	\begin{equation} \label{eqnlae2}
		\mathbf{Y}
		\begin{pmatrix}
			\tau^*  s_1 \\
			\vdots \\
			\tau^*  s_n 
		\end{pmatrix}
		= \begin{pmatrix}
			\tau^*  \\
			\vdots \\
			\tau^*
		\end{pmatrix}
		= \Vert P_0 \Vert^2
		\begin{pmatrix}
			1 \\
			\vdots \\
			1
		\end{pmatrix}
	\end{equation}
	where we recall (\ref{eqtst}).
	
	If $P_0$ is in the interior of $\mathrm{Conv} (\CF )$, (\ref{eqnlae1}) and (\ref{eqnlae2}) immediately imply that (\ref{eqnikyyt}) admits a solution $\alpha_1 , \dots , \alpha_m >0 $ by setting $\alpha_i := \Vert P_0 \Vert^{-2} \beta_i >0$ for $i=1, \dots , m$, which indeed satisfies $\sum_{i=1}^m \alpha_i = \Vert P_0 \Vert^{-2}$. Conversely, if (\ref{eqnikyyt}) admits a solution $\alpha_1 , \dots , \alpha_m >0$ satisfying $\sum_{i=1}^m \alpha_i = \Vert P_0 \Vert^{-2}$, then (\ref{eqnlae2}) implies
	\begin{equation} \label{eqnlae3}
		\mathbf{Y} \left(
		\mathbf{Y}^t
		\begin{pmatrix}
			\alpha_1  \\
			\vdots \\
			\alpha_m 
		\end{pmatrix}
		- \Vert P_0 \Vert^{-2}
		\begin{pmatrix}
			\tau^*  s_1 \\
			\vdots \\
			\tau^*  s_n 
		\end{pmatrix}
		\right)
		=
		\Vert P_0 \Vert^{-2} \mathbf{Y} \left(
		\mathbf{Y}^t
		\begin{pmatrix}
			\alpha_1 \Vert P_0 \Vert^{2} \\
			\vdots \\
			\alpha_m \Vert P_0 \Vert^{2}
		\end{pmatrix}
		- P_0
		\right) =0.
	\end{equation}
	We set
	\begin{equation*}
	w_0 := \mathbf{Y}^t \begin{pmatrix} \alpha_1 \Vert P_0 \Vert^{2} \\ \vdots \\ \alpha_m \Vert P_0 \Vert^{2} \end{pmatrix},
	\end{equation*}
	which is a strict convex combination of elements in $\CF$ and hence is in the interior of $\mathrm{Conv} (\CF )$. The equation (\ref{eqnlae3}) implies that $w_0 - P_0$ is orthogonal to any element in $\CF$, which in particular implies $\langle w_0 , w_0 - P_0 \rangle = 0$ since $w_0 \in \mathrm{Conv} (\CF )$. We also have $\langle  P_0 , w_0 - P_0 \rangle =0$ by (\ref{eqdftsot}) and (\ref{eqdftsot2}), together with $w_0 \in \mathrm{Conv} (\CF )$. Thus we get
	\begin{equation*}
		\Vert w_0 - P_0 \Vert^2 = \langle w_0 - P_0 , w_0 - P_0 \rangle = \langle w_0 , w_0 - P_0 \rangle - \langle  P_0 , w_0 - P_0 \rangle =0,
	\end{equation*}
	which implies $w_0 = P_0$, and hence $P_0$ is in the interior of $\mathrm{Conv} (\CF )$.
\end{proof}

\begin{remark} \label{rmnkscal}
	Nikolayevsky \cite{Nik08ded,Nik11} does not seem to impose the condition $\sum_{i=1}^m \alpha_i = \Vert P_0 \Vert^{-2}$ for the solution to (\ref{eqnikyyt}), at least explicitly, and it could be redundant. We added this condition in the statement of Theorem \ref{thnik} since in any case it is always satisfied when $P_0$ is in the interior of $\mathrm{Conv} (\CF )$, as we can see from the above proof.
\end{remark}

\subsection{Modified Taketomi--Tamaru conjecture} \label{scmttc}

Taketomi--Tamaru \cite{TakTam18} conjectured the following, and exhibited examples for which it holds.

\begin{conjecture} (Taketomi--Tamaru \cite[Conjecture 1.2]{TakTam18})
	Let $\Fl$ be an $n$-dimensional Lie algebra, and $L$ be the connected and simply connected Lie group with the Lie algebra $\Fl$. If the (right) action of $\rl^{\times} \mathrm{Aut} (\Fl )$ on $O(n) \backslash GL_n (\rl)$ is not transitive, and all $\rl^{\times} \mathrm{Aut} (\Fl )$-orbits are congruent to each other with respect to $GL_n (\rl)$, then $L$ does not admit left-invariant Ricci solitons.
\end{conjecture}

Unfortunately the claim does not hold in the stated form (at least when we consider the identity component of $\mathrm{Aut} (\Fl )$) since there is a counterexample by Jablonski \cite{Jab18}. We propose and prove a significantly modified version of the above conjecture when $L$ is a nilpotent Lie group, as an application of the variational formulation discussed in \S \ref{sckne}.

\begin{theorem} \label{ppttcjw}
	Let $(N , \textsl{g}_N)$ be an $n$-dimensional homogeneous nilmanifold with the metric Lie algebra $(\Fn , \mu , \langle , \rangle )$, compatible with the pre-Einstein derivation, and $G$, $K$, $\Fp$ be as defined in Definitions \ref{dflagk} and \ref{dfgpgk}. Suppose that the following hold:
	\begin{enumerate}
		\item the (right) action of $\mathrm{Aut} (\mu ) \cap G$ on the right coset space $Y = K \backslash G$ is not transitive,
		\item all $\mathrm{Aut} ( \mu ) \cap G$-orbits are congruent to each other with respect to $G$, i.e.~for any $p_1, p_2 \in Y$, there exists $g \in G$ such that $p_1 \cdot ( \mathrm{Aut} ( \mu ) \cap G ) = (p_2 \cdot (\mathrm{Aut} ( \mu ) \cap G)) \cdot g$,
		\item for any $g \in G \setminus  \mathrm{Aut} (\mu ) $, there exists a non-zero element in $\Fp \cap ( \mathrm{Der} (\mu) \setminus \mathrm{Der} (\rho (g) \cdot \mu))$.
	\end{enumerate}
	Then the Kempf--Ness energy is not bounded from below; in particular, $N$ cannot admit left-invariant Ricci solitons.
\end{theorem}

\begin{proof}
	Suppose for contradiction that $E_{\mu} : Y \to \rl$ is bounded from below; note that $E_{\mu}$ is not a constant function since the action of $\mathrm{Aut} (\mu ) \cap G$ on $Y$ is not transitive by the first hypothesis (otherwise it contradicts the strict convexity in Lemma \ref{lmkncvx}).
	
	Pick any $g_* \in G \setminus \mathrm{Aut} ( \mu)$, and take $p_1 := K \cdot g_*$ and $p_2 := K \cdot I_{\Fn}$. Then $\rho ( p_1 \cdot ( \mathrm{Aut} ( \mu ) \cap G )) \cdot \mu = \rho (K ) \cdot (\rho (g_*) \cdot \mu)$. By the second hypothesis there exists $g \in G$ such that $p_1 \cdot ( \mathrm{Aut} ( \mu ) \cap G) = (p_2 \cdot ( \mathrm{Aut} ( \mu ) \cap G)) \cdot g$, which implies
	\begin{equation*} 
		\rho (K) \cdot (\rho(g_*) \cdot \mu) = \rho ( K \cdot ( \mathrm{Aut} ( \mu) \cap G )) \cdot ( \rho (g) \cdot \mu ).
	\end{equation*}
	In particular, this implies that for any $h \in \mathrm{Aut} ( \mu) \cap G$ and the point $K \cdot hg \in Y$ which corresponds to $hg$, we have
	\begin{equation*}
		E_{\mu} (K \cdot hg) = E_{\mu} (K \cdot g_*) 
	\end{equation*}
	which is constant in $h$. Thus 
	\begin{equation} \label{eqknhgc}
		E_{\rho (g) \cdot \mu} (K \cdot h) = \log \Vert \rho (h) \cdot \rho (g) \cdot \mu \Vert^2 =  \log \Vert \rho (hg) \cdot \mu \Vert^2 =  E_{\mu} (K \cdot hg) = E_{\mu} (K \cdot g_*)
	\end{equation}
	holds for all $h \in  \mathrm{Aut} ( \mu) \cap G$.

	We now take a geodesic ray $\gamma_{\lambda} (t) := K \cdot \exp (\lambda t /2 ) \subset \mathrm{Aut} ( \mu ) \cap G$, generated by a non-zero element $\lambda \in \Fp \cap (\mathrm{Der} (\mu) ) \setminus \mathrm{Der} (\rho (g) \cdot \mu)$, which exists by the third hypothesis. Differentiating the equation (\ref{eqknhgc}) along $\gamma_{\lambda} (t)$, we find
	\begin{equation*}
		\frac{d^2}{dt^2} E_{\rho (g) \cdot \mu} ( \gamma_{\lambda} (t)) = 0
	\end{equation*}
	for all $t \ge 0$, but this forces $\lambda \in \mathrm{Der} (\rho(g) \cdot \mu )$ by the strict convexity (Lemma \ref{lmkncvx}) since we assumed that $E_{\mu}$, and hence $E_{\rho (g) \cdot \mu}$, is bounded from below over $Y$ (see also Lemma \ref{ppsvspps}). Thus we find $\lambda \in \mathrm{Der} (\mu ) \cap \mathrm{Der} (\rho(g) \cdot \mu )$, which contradicts the assumption $\lambda \in \mathrm{Der} (\mu)  \setminus \mathrm{Der} (\rho (g) \cdot \mu)$.
\end{proof}

Theorem \ref{ppttcjw} above is weaker and more complicated than the original conjecture, partly because of the appearance of more complicated groups $G$ and $K$ (and also $(\Fn , \mu)$ is assumed to be nilpotent), but the author believes that addition of the third condition is the most significant as it is not met by Jablonski's counterexample. In the notation of his example \cite[Theorem 2.2]{Jab18}, we can easily compute
\begin{equation*}
	\CS = \rl \left( I_{\Fn} - \frac{2 \dim \Fz + \dim \Fv}{4 \dim \Fz + \dim \Fv} \begin{bmatrix}  I_{\Fv} & 0 \\ 0 & 2 I_{\Fz} \end{bmatrix} \right)  = \rl \begin{bmatrix} 2\dim \Fz I_{\Fv} & 0 \\ 0 & - \dim \Fv I_{\Fz} \end{bmatrix}.
\end{equation*}
When we take
\begin{equation*}
	g := \begin{bmatrix} D_1 & 0 \\ 0 & D_2 \end{bmatrix} ,
\end{equation*}
where $D_1$ (resp.~$D_2$) is any diagonal matrix of unit determinant acting on $\Fv$ (resp.~$\Fz$), we find $g \in G \setminus \mathrm{Aut} (\mu )$ if $D_1$ and $D_2$ are chosen generically, and that $g$ fixes the $(1, 2)$-derivation $D$ in \cite[Theorem 2.2]{Jab18}. In particular, there exist no non-zero self-adjoint elements in $\mathrm{Der} (\mu) \setminus \mathrm{Der} (\rho (g) \cdot \mu) = \mathrm{Der} (\mu) \setminus g \mathrm{Der} ( \mu) g^{-1}$, and hence $\Fp \cap ( \mathrm{Der} (\mu) \setminus \mathrm{Der} (\rho (g) \cdot \mu) ) =0$ for Jablonski's counterexample. Noting that his example has a ``very small'' derivation algebra, it seems natural to wonder if the third condition in Theorem \ref{ppttcjw} is satisfied generically when the derivation algebra is relatively large, but we do not pursue any further results in this paper.

\bibliography{nilsoliton.bib}

\begin{bibdiv}
\begin{biblist}

\bib{BilWin23}{article}{
      author={Biliotti, Leonardo},
      author={Windare, Oluwagbenga~Joshua},
       title={A {H}ilbert-{M}umford criterion for polystability for actions of
  real reductive lie groups},
        date={2023},
     journal={arXiv preprint arXiv:2309.01138},
}

\bib{BilWin23p}{article}{
      author={Biliotti, Leonardo},
      author={Windare, Oluwagbenga~Joshua},
       title={Stability, analytic stability for real reductive {L}ie groups},
        date={2023},
        ISSN={1050-6926},
     journal={J. Geom. Anal.},
      volume={33},
      number={3},
       pages={Paper No. 92, 31},
         url={https://doi.org/10.1007/s12220-022-01146-0},
      review={\MR{4531069}},
}

\bib{BL22}{article}{
      author={B\"{o}hm, Christoph},
      author={Lafuente, Ramiro~A.},
       title={Homogeneous {E}instein metrics on {E}uclidean spaces are
  {E}instein solvmanifolds},
        date={2022},
        ISSN={1465-3060},
     journal={Geom. Topol.},
      volume={26},
      number={2},
       pages={899\ndash 936},
         url={https://doi.org/10.2140/gt.2022.26.899},
      review={\MR{4444271}},
}

\bib{BL23}{article}{
      author={B\"{o}hm, Christoph},
      author={Lafuente, Ramiro~A.},
       title={Non-compact {E}instein manifolds with symmetry},
        date={2023},
        ISSN={0894-0347},
     journal={J. Amer. Math. Soc.},
      volume={36},
      number={3},
       pages={591\ndash 651},
         url={https://doi.org/10.1090/jams/1022},
      review={\MR{4583772}},
}

\bib{BHJ1}{article}{
      author={Boucksom, S\'{e}bastien},
      author={Hisamoto, Tomoyuki},
      author={Jonsson, Mattias},
       title={Uniform {K}-stability, {D}uistermaat-{H}eckman measures and
  singularities of pairs},
        date={2017},
        ISSN={0373-0956},
     journal={Ann. Inst. Fourier (Grenoble)},
      volume={67},
      number={2},
       pages={743\ndash 841},
         url={http://aif.cedram.org/item?id=AIF_2017__67_2_743_0},
      review={\MR{3669511}},
}

\bib{Dol}{book}{
      author={Dolgachev, Igor},
       title={Lectures on invariant theory},
      series={London Mathematical Society Lecture Note Series},
   publisher={Cambridge University Press, Cambridge},
        date={2003},
      volume={296},
        ISBN={0-521-52548-9},
         url={https://doi.org/10.1017/CBO9780511615436},
      review={\MR{2004511}},
}

\bib{EbeJab09}{incollection}{
      author={Eberlein, Patrick},
      author={Jablonski, Michael},
       title={Closed orbits of semisimple group actions and the real
  {H}ilbert-{M}umford function},
        date={2009},
   booktitle={New developments in {L}ie theory and geometry},
      series={Contemp. Math.},
      volume={491},
   publisher={Amer. Math. Soc., Providence, RI},
       pages={283\ndash 321},
         url={https://doi.org/10.1090/conm/491/09620},
      review={\MR{2537062}},
}

\bib{FerCul14}{article}{
      author={Fern\'{a}ndez-Culma, Edison~Alberto},
       title={Classification of nilsoliton metrics in dimension seven},
        date={2014},
        ISSN={0393-0440},
     journal={J. Geom. Phys.},
      volume={86},
       pages={164\ndash 179},
         url={https://doi.org/10.1016/j.geomphys.2014.07.032},
      review={\MR{3282320}},
}

\bib{FutMab95}{article}{
      author={Futaki, Akito},
      author={Mabuchi, Toshiki},
       title={Bilinear forms and extremal {K}\"{a}hler vector fields associated
  with {K}\"{a}hler classes},
        date={1995},
        ISSN={0025-5831},
     journal={Math. Ann.},
      volume={301},
      number={2},
       pages={199\ndash 210},
         url={https://doi.org/10.1007/BF01446626},
      review={\MR{1314584}},
}

\bib{HasHaya}{incollection}{
      author={Hashimoto, Yoshinori},
       title={Balanced metrics for extremal {K}\"{a}hler metrics and {F}ano
  manifolds},
        date={[2024] \copyright 2024},
   booktitle={The {B}ergman kernel and related topics},
      series={Springer Proc. Math. Stat.},
      volume={447},
   publisher={Springer, Singapore},
       pages={169\ndash 188},
         url={https://doi.org/10.1007/978-981-99-9506-6_5},
      review={\MR{4731752}},
}

\bib{Heb98}{article}{
      author={Heber, Jens},
       title={Noncompact homogeneous {E}instein spaces},
        date={1998},
        ISSN={0020-9910},
     journal={Invent. Math.},
      volume={133},
      number={2},
       pages={279\ndash 352},
         url={https://doi.org/10.1007/s002220050247},
      review={\MR{1632782}},
}

\bib{Hor07}{book}{
      author={H\"{o}rmander, Lars},
       title={Notions of convexity},
      series={Modern Birkh\"{a}user Classics},
   publisher={Birkh\"{a}user Boston, Inc., Boston, MA},
        date={2007},
        ISBN={978-0-8176-4584-7; 0-8176-4584-5},
        note={Reprint of the 1994 edition [of MR1301332]},
      review={\MR{2311920}},
}

\bib{Jab10}{article}{
      author={Jablonski, Michael},
       title={Detecting orbits along subvarieties via the moment map},
        date={2010},
        ISSN={1867-5778},
     journal={M\"{u}nster J. Math.},
      volume={3},
       pages={67\ndash 88},
      review={\MR{2775356}},
}

\bib{Jab11}{article}{
      author={Jablonski, Michael},
       title={Concerning the existence of {E}instein and {R}icci soliton
  metrics on solvable {L}ie groups},
        date={2011},
        ISSN={1465-3060},
     journal={Geom. Topol.},
      volume={15},
      number={2},
       pages={735\ndash 764},
         url={https://doi.org/10.2140/gt.2011.15.735},
      review={\MR{2800365}},
}

\bib{Jab14}{article}{
      author={Jablonski, Michael},
       title={Homogeneous {R}icci solitons are algebraic},
        date={2014},
        ISSN={1465-3060},
     journal={Geom. Topol.},
      volume={18},
      number={4},
       pages={2477\ndash 2486},
         url={https://doi.org/10.2140/gt.2014.18.2477},
      review={\MR{3268781}},
}

\bib{Jab23}{article}{
      author={Jablonski, Michael},
       title={Homogeneous {E}instein manifolds},
        date={2023},
        ISSN={0041-6932},
     journal={Rev. Un. Mat. Argentina},
      volume={64},
      number={2},
       pages={461\ndash 485},
      review={\MR{4771217}},
}

\bib{Jab18}{article}{
      author={Jablonski, Michael},
       title={Concerning a conjecture of {T}aketomi-{T}amaru},
        date={2025},
     journal={arXiv preprint arXiv:1810.08173v2},
}

\bib{Jen69}{article}{
      author={Jensen, Gary~R.},
       title={Homogeneous {E}instein spaces of dimension four},
        date={1969},
        ISSN={0022-040X},
     journal={J. Differential Geometry},
      volume={3},
       pages={309\ndash 349},
         url={http://projecteuclid.org/euclid.jdg/1214429056},
      review={\MR{261487}},
}

\bib{KemNes78}{inproceedings}{
      author={Kempf, George~R.},
      author={Ness, Linda~A.},
       title={The length of vectors in representation spaces},
        date={1979},
   booktitle={Algebraic geometry ({P}roc. {S}ummer {M}eeting, {U}niv.
  {C}openhagen, {C}openhagen, 1978)},
      series={Lecture Notes in Math.},
      volume={732},
   publisher={Springer, Berlin},
       pages={233\ndash 243},
      review={\MR{555701}},
}

\bib{Kirw}{book}{
      author={Kirwan, Frances~Clare},
       title={Cohomology of quotients in symplectic and algebraic geometry},
      series={Mathematical Notes},
   publisher={Princeton University Press, Princeton, NJ},
        date={1984},
      volume={31},
        ISBN={0-691-08370-3},
         url={https://doi.org/10.2307/j.ctv10vm2m8},
      review={\MR{766741}},
}

\bib{KN2}{book}{
      author={Kobayashi, Shoshichi},
      author={Nomizu, Katsumi},
       title={Foundations of differential geometry. {V}ol. {II}},
      series={Wiley Classics Library},
   publisher={John Wiley \& Sons, Inc., New York},
        date={1996},
        ISBN={0-471-15732-5},
        note={Reprint of the 1969 original, A Wiley-Interscience Publication},
      review={\MR{1393941}},
}

\bib{Lau01}{article}{
      author={Lauret, Jorge},
       title={Ricci soliton homogeneous nilmanifolds},
        date={2001},
        ISSN={0025-5831},
     journal={Math. Ann.},
      volume={319},
      number={4},
       pages={715\ndash 733},
         url={https://doi.org/10.1007/PL00004456},
      review={\MR{1825405}},
}

\bib{Lau01b}{article}{
      author={Lauret, Jorge},
       title={Standard {E}instein solvmanifolds as critical points},
        date={2001},
        ISSN={0033-5606},
     journal={Q. J. Math.},
      volume={52},
      number={4},
       pages={463\ndash 470},
         url={https://doi.org/10.1093/qjmath/52.4.463},
      review={\MR{1874492}},
}

\bib{Lau02}{article}{
      author={Lauret, Jorge},
       title={Finding {E}instein solvmanifolds by a variational method},
        date={2002},
        ISSN={0025-5874},
     journal={Math. Z.},
      volume={241},
      number={1},
       pages={83\ndash 99},
         url={https://doi.org/10.1007/s002090100407},
      review={\MR{1930986}},
}

\bib{Lau03}{article}{
      author={Lauret, Jorge},
       title={On the moment map for the variety of {L}ie algebras},
        date={2003},
        ISSN={0022-1236},
     journal={J. Funct. Anal.},
      volume={202},
      number={2},
       pages={392\ndash 423},
         url={https://doi.org/10.1016/S0022-1236(02)00108-8},
      review={\MR{1990531}},
}

\bib{Lau06}{article}{
      author={Lauret, Jorge},
       title={A canonical compatible metric for geometric structures on
  nilmanifolds},
        date={2006},
        ISSN={0232-704X},
     journal={Ann. Global Anal. Geom.},
      volume={30},
      number={2},
       pages={107\ndash 138},
         url={https://doi.org/10.1007/s10455-006-9015-y},
      review={\MR{2234091}},
}

\bib{Lau09}{incollection}{
      author={Lauret, Jorge},
       title={Einstein solvmanifolds and nilsolitons},
        date={2009},
   booktitle={New developments in {L}ie theory and geometry},
      series={Contemp. Math.},
      volume={491},
   publisher={Amer. Math. Soc., Providence, RI},
       pages={1\ndash 35},
         url={https://doi.org/10.1090/conm/491/09607},
      review={\MR{2537049}},
}

\bib{Lau10}{article}{
      author={Lauret, Jorge},
       title={Einstein solvmanifolds are standard},
        date={2010},
        ISSN={0003-486X},
     journal={Ann. of Math. (2)},
      volume={172},
      number={3},
       pages={1859\ndash 1877},
         url={https://doi.org/10.4007/annals.2010.172.1859},
      review={\MR{2726101}},
}

\bib{Lau11}{article}{
      author={Lauret, Jorge},
       title={Ricci soliton solvmanifolds},
        date={2011},
        ISSN={0075-4102},
     journal={J. Reine Angew. Math.},
      volume={650},
       pages={1\ndash 21},
         url={https://doi.org/10.1515/CRELLE.2011.001},
      review={\MR{2770554}},
}

\bib{LauWil11}{article}{
      author={Lauret, Jorge},
      author={Will, Cynthia},
       title={Einstein solvmanifolds: existence and non-existence questions},
        date={2011},
        ISSN={0025-5831},
     journal={Math. Ann.},
      volume={350},
      number={1},
       pages={199\ndash 225},
         url={https://doi.org/10.1007/s00208-010-0552-0},
      review={\MR{2785768}},
}

\bib{Mar01}{article}{
      author={Marian, Alina},
       title={On the real moment map},
        date={2001},
        ISSN={1073-2780},
     journal={Math. Res. Lett.},
      volume={8},
      number={5-6},
       pages={779\ndash 788},
         url={https://doi.org/10.4310/MRL.2001.v8.n6.a8},
      review={\MR{1879820}},
}

\bib{Mil76}{article}{
      author={Milnor, John},
       title={Curvatures of left invariant metrics on {L}ie groups},
        date={1976},
        ISSN={0001-8708},
     journal={Advances in Math.},
      volume={21},
      number={3},
       pages={293\ndash 329},
         url={https://doi.org/10.1016/S0001-8708(76)80002-3},
      review={\MR{425012}},
}

\bib{MFK94}{book}{
      author={Mumford, D.},
      author={Fogarty, J.},
      author={Kirwan, F.},
       title={Geometric invariant theory},
     edition={Third},
      series={Ergebnisse der Mathematik und ihrer Grenzgebiete (2) [Results in
  Mathematics and Related Areas (2)]},
   publisher={Springer-Verlag, Berlin},
        date={1994},
      volume={34},
        ISBN={3-540-56963-4},
      review={\MR{1304906}},
}

\bib{Nik08ded}{article}{
      author={Nikolayevsky, Yuri},
       title={Einstein solvmanifolds with a simple {E}instein derivation},
        date={2008},
        ISSN={0046-5755},
     journal={Geom. Dedicata},
      volume={135},
       pages={87\ndash 102},
         url={https://doi.org/10.1007/s10711-008-9264-y},
      review={\MR{2413331}},
}

\bib{Nik11}{article}{
      author={Nikolayevsky, Yuri},
       title={Einstein solvmanifolds and the pre-{E}instein derivation},
        date={2011},
        ISSN={0002-9947},
     journal={Trans. Amer. Math. Soc.},
      volume={363},
      number={8},
       pages={3935\ndash 3958},
         url={https://doi.org/10.1090/S0002-9947-2011-05045-2},
      review={\MR{2792974}},
}

\bib{Pay10}{article}{
      author={Payne, Tracy~L.},
       title={The existence of soliton metrics for nilpotent {L}ie groups},
        date={2010},
        ISSN={0046-5755},
     journal={Geom. Dedicata},
      volume={145},
       pages={71\ndash 88},
         url={https://doi.org/10.1007/s10711-009-9404-z},
      review={\MR{2600946}},
}

\bib{RicSlo90}{article}{
      author={Richardson, R.~W.},
      author={Slodowy, P.~J.},
       title={Minimum vectors for real reductive algebraic groups},
        date={1990},
        ISSN={0024-6107},
     journal={J. London Math. Soc. (2)},
      volume={42},
      number={3},
       pages={409\ndash 429},
         url={https://doi.org/10.1112/jlms/s2-42.3.409},
      review={\MR{1087217}},
}

\bib{Szethesis}{thesis}{
      author={Sz\'{e}kelyhidi, G\'{a}bor},
       title={{E}xtremal metrics and {K}-stability},
        type={Ph.D. Thesis},
        date={2006},
}

\bib{Sze07}{article}{
      author={Sz\'{e}kelyhidi, G\'{a}bor},
       title={Extremal metrics and {$K$}-stability},
        date={2007},
        ISSN={0024-6093},
     journal={Bull. Lond. Math. Soc.},
      volume={39},
      number={1},
       pages={76\ndash 84},
         url={https://doi.org/10.1112/blms/bdl015},
      review={\MR{2303522}},
}

\bib{TakTam18}{article}{
      author={Taketomi, Y.},
      author={Tamaru, H.},
       title={On the nonexistence of left-invariant {R}icci solitons---a
  conjecture and examples},
        date={2018},
        ISSN={1083-4362},
     journal={Transform. Groups},
      volume={23},
      number={1},
       pages={257\ndash 270},
         url={https://doi.org/10.1007/s00031-017-9439-4},
      review={\MR{3763948}},
}

\bib{TiaZhu02}{article}{
      author={Tian, Gang},
      author={Zhu, Xiaohua},
       title={A new holomorphic invariant and uniqueness of
  {K}\"{a}hler-{R}icci solitons},
        date={2002},
        ISSN={0010-2571},
     journal={Comment. Math. Helv.},
      volume={77},
      number={2},
       pages={297\ndash 325},
         url={https://doi.org/10.1007/s00014-002-8341-3},
      review={\MR{1915043}},
}

\bib{Wil03}{article}{
      author={Will, Cynthia},
       title={Rank-one {E}instein solvmanifolds of dimension 7},
        date={2003},
        ISSN={0926-2245},
     journal={Differential Geom. Appl.},
      volume={19},
      number={3},
       pages={307\ndash 318},
         url={https://doi.org/10.1016/S0926-2245(03)00037-8},
      review={\MR{2013098}},
}

\bib{Wil82}{article}{
      author={Wilson, Edward~N.},
       title={Isometry groups on homogeneous nilmanifolds},
        date={1982},
        ISSN={0046-5755},
     journal={Geom. Dedicata},
      volume={12},
      number={3},
       pages={337\ndash 346},
         url={https://doi.org/10.1007/BF00147318},
      review={\MR{661539}},
}

\end{biblist}
\end{bibdiv}

\end{document}